\documentclass{article}

\usepackage[american]{babel}
\usepackage{amsmath,amssymb,theorem, graphics, graphicx, wrapfig,relsize}  
\usepackage{color}

\usepackage{url} 

\newcommand{\blind}{1}

\newcommand{\qed}{\ \hglue 0pt plus 1filll $\Box$}

\newcommand{\fK}{\mathfrak{K}}
\newcommand{\fb}{\mathfrak{b}}
\newcommand{\bE}{\mathbb{E}}
\newcommand{\bN}{\mathbb{N}}
\newcommand{\bP}{\mathbb{P}}
\newcommand{\bR}{\mathbb{R}}

\newcommand{\cC}{\mathcal{C}}
\newcommand{\cD}{\mathcal{D}}
\newcommand{\cF}{\mathcal{F}}
\newcommand{\cN}{\mathcal{N}}

\newcommand{\ep}{\varepsilon}

\newcommand{\var}{\text{Var}}

\newcommand{\set}[1]{\left\{#1 \right\}}
\newcommand{\bra}[1]{\left(#1 \right)}
\newcommand{\brac}[1]{\left[#1 \right]}
\newcommand{\ubra}[1]{\lfloor #1/H \rfloor\cdot H}
\newcommand{\abs}[1]{\left|#1\right|}

\newcommand{\ot}[1]{\overline{\triangle #1}}
\newcommand{\con}{ \ \Big| \ }

\theoremstyle{change}
\newtheorem{thm}{Theorem}[section]
\newtheorem{prop}[thm]{Proposition}
\newtheorem{exa}[thm]{Example}
\newtheorem{rem}[thm]{Remark}
\newtheorem{assu}[thm]{Assumption}

\usepackage{titlesec}
\titleformat*{\section}{\large\bfseries}
\titleformat*{\subsection}{\bfseries}
\titlelabel{\thesection. \ }

\begin{document}







\if1\blind
{
  \title{\bf Volatility Decomposition and Estimation in Time-Changed Price Models}
  \author{Rainer Dahlhaus\thanks{
    This work was supported by the Deutsche Foschungsgemeinschaft, RTG 1953. The authors thank Per Mykland and Lan Zhang for very helpful comments.}\hspace{.2cm}\\
    Institute of Applied Mathematics, Heidelberg University\\
   	Im Neuenheimer Feld 205, 69120\\
   	Heidelberg, Germany\\
    E-mail: dahlhaus@statlab.uni-heidelberg.de\\
    and \\
    Sophon Tunyavetchakit  \\
    Institute of Applied Mathematics, Heidelberg University\\
    Im Neuenheimer Feld 205, 69120\\
   	Heidelberg, Germany\\
    E-mail: s.tunyavetchakit@uni-heidelberg.de.
    }
  \maketitle
} \fi

\if0\blind
{
  \bigskip
  \bigskip
  \bigskip
  \begin{center}
    {\LARGE\bf Volatility Decomposition and Estimation in Time-Changed Price Models}
\end{center}
  \medskip
} \fi

\bigskip
\begin{abstract}
The usage of a spot volatility estimate based on a volatility decomposition in a time-changed price-model according to the trading times is investigated. In this model clock-time volatility splits up into the product of tick-time volatility and trading intensity, which both can be estimated from data and contain valuable information. By inspecting these two curves individually we gain more insight into the cause and structure of volatility. Several examples are provided where the tick-time volatility curve is much smoother than the clock-time volatility curve meaning that the major part of fluctuations in clock-time volatility is due to fluctuations of the trading intensity. Since microstructure noise only influences the estimation of the (smooth) tick-time volatility curve, the findings lead to an improved pre-averaging estimator of spot volatility. This is reflected by a better rate of convergence of the estimator. The asymptotic properties of the estimators are derived by an infill asymptotic approach.
\end{abstract}

\noindent%
\vfill

\newpage

\section{INTRODUCTION} \label{ch:introduction}

The estimation of volatility for high-frequency data under microstructure noise has been extensively studied during recent years - see A\"it-Sahalia and Jacod~(2014) for an overview. The majority of this work has been carried out in the framework of diffusion models. In this paper we focus on the estimation of spot volatility, and, contrary to the majority of previous research,  on a time-changed price-model based on trading times. In this model volatility splits up into the product of two identifiable curves, namely tick-time volatility and trading intensity. The main methodological and theoretical contributions of this paper are the introduction and the theoretical investigation of a volatility estimate based on this volatility decomposition, and the proof that this estimator can outperform the classical (diffusion model based) estimators in terms of the rate of convergence.

For the estimation of tick-time volatility under microstructure noise we have to adapt an estimator from diffusion models to our situation which can cope successfully with microstructure noise in a high-frequency situation. Many noise-robust estimators have been introduced in the literature. In Zhang et al.~(2005) the combination of two different timescales is used to construct a consistent estimator for the integrated volatility. This idea is extended later to the multi-timescale estimator that archives the optimal rate of convergence $n^{-1/4}$ - cf.~Zhang~(2006). Barndorff-Nielsen et~al.~(2008) suggest a flat-top kernel-type estimator, called realized kernel, which combines different lags of autocovariances to eliminate the effect of microstructure noise. Recently, an estimator presented by Rei\ss \ (2011) and Bibinger and Rei\ss \ (2014) has received attention since their estimator is asymptotically efficient. In this work, we apply the pre-averaging technique, which was introduced by Podolskij and Vetter~(2009) and later extended by Jacod et al.~(2009), in order to construct a noise-robust estimate for the tick-time volatility in our time-changed model. However, it is also possible to adapt most of the other methods to the model of this paper. An estimator based on particle filtering in a nonlinear microstructure noise model has been discussed in Dahlhaus and Neddermeyer~(2013).

Time-changed price-models were first investigated in Clark~(1973) in connection with finance. In his work, the volume of trades is suggested to be a subordinator of a Brownian motion in order to recover the normality of the distribution of cotton future prices. Afterwards a relationship between asset returns, price fluctuation, and market activities measured by trading volume and numbers of transactions is extensively discussed. An\'e and Geman~(2000) conclude that the number of trades explains the volatility change better than their volume, so they recover the normality of asset returns through this stochastic time change in high-frequency data; see also Jones et al.~(1994), Plerou et al.~(2001) and Gabaix et al.~(2003) for more detailed discussions of this correlation. Due to various mathematical tools, the time-changed Brownian motion is attractive and tractable to study arbitrage-free asset returns, which are shown to be semimartingales (see e.g.~Delbaen and Schachermayer~1994). Indeed, having a class of time-changed Brownian motion is satisfactory since it is as large as a class of semimartingale; see Monroe~(1978). In recent years, other time change models have been extensively studied, especially a time-changed L\'evy process which allows for a more complex structure in the price models coping with some stylized-effect emerging in the real market; for details refer to Carr et al.~(2003) and Carr and Wu~(2004), and to Belomestny~(2011) for a statistical treatment of this kind of models.

The article is organized as follows. Section~\ref{ch:voladecomposition} contains an introduction to the model, the volatility decomposition, and the estimates with a discussion of the implications for applications. Here, the volatility decomposition is proven under a general setup. In Section~\ref{ch:InfillAsymptotics} we investigate the asymptotic properties of the estimates by means of an infill asymptotic approach constructed by time-rescaling. The asymptotic results are compared in Table~\ref{table:Comparison_PAVG_Volatility} which shows the advantage of using the volatility decomposition. In particular the results show that the rate of convergence of classical estimators can be outperformed within the model of this paper. Section~\ref{ch:conclusion} contains some concluding remarks. The details of the data analysis and the proofs can be found in the Appendix. If not otherwise stated, all equalities and inequalities of random expressions are in an almost sure sense.

\section{THE VOLATILITY DECOMPOSITION} \label{ch:voladecomposition}

As motivated in the introduction we use instead of the classical semimartingale model a diffusion model subordinated by transaction time, for example the time-changed Brownian motion
\begin{equation} \label{TransactionTime-Model}
dX_t = \sigma_t \, dW_{N_t} \quad \text{ for } t \in [0,T]
\end{equation}
with $N_t$ being a point process with intensity $\lambda_t$ reflecting the accumulated number of transactions up to time $t$. In it's simplest form $\sigma_t$ and $\lambda_t$ are deterministic and $W(\cdot)$ and $N_{\cdot}$ are independent. In a more general model $\sigma_t$ and $\lambda_t$ are stochastic processes depending on the past of $X_t$ and $N_t$, the independence of $W(\cdot)$ and $N_{\cdot}$ may be replaced by some martingale-structure, and $W(\cdot)$ may be non-Gaussian. We include microstructure noise into our considerations - for example the asymptotic properties of our estimates are derived under the assumption of additive i.i.d. noise
\begin{equation} \label{}
Y_{t_i} = X_{t_i} + \ep_i  \quad \text{ for } i=1,...,N_T,
\end{equation}
where $t_i$ are the trading times.

The focus of this paper is the estimation of spot volatility for financial transactions which is \underline{not} the function $\sigma^2_t$ from the above model. Even more the meaning of $\sigma^2_t$ is different from the meaning of $\sigma^2_{\!clock} (t)$ in the classical ``clock-time'' diffusion model (say $dX_t = \sigma_{\!clock} (t) \, dW_t$). We therefore start with a model-independent definition of spot volatility and clarify the relation to $\sigma_t$ in the different models. Let $(\cF_{t})_{t \geq 0}$ be an increasing sequence of $\sigma$-algebras; roughly speaking, it represents the information available up to and including time $t$. We define
\begin{equation*} \label{}
\rm{vola}^{2}_{\,t} := \lim_{\Delta t \rightarrow 0} \frac
{\bE \big[ (\mathbf{X}(t+\Delta t) - \mathbf{X}(t))^{2} \big| \cF_{t} \big]} {\Delta t}\,.
\end{equation*}
In the classical diffusion model, $\rm{vola}^{2}_{\,t}$ does not depend on the point process and we have under the assumption that $\sigma_{\!clock} (t)$ is a right-continuous process with left-limits adapted to $\cF_{t}$:  $\rm{vola}^{2}_{\,t} = \sigma_{clock}^{2}(t)$. Therefore, we use $\sigma_{clock}(t)$ in this paper as a synonym for $\rm{vola}_{\,t}$, i.e.~we define
\begin{equation*} \label{}
\sigma_{clock}^{2}(t) := \rm{vola}^{2}_{\,t}.
\end{equation*}
In the transaction-time model of this paper we prove below that $\sigma_{clock}^{2}(t) = \sigma_t^{2} \cdot \lambda_t$. We first set down the assumptions for this result. To understand our assumptions, note that in (\ref{TransactionTime-Model}) we do not use the whole process $W(\cdot)$ but only the increments $U_i := W (N_{t_i}) - W(N_{t_{i-1}})$ which we now assume to be a martingale difference sequence.

\begin{assu} \label{assu_2_1}
The $X_{t_i}$ at observation times $t_i$ follow the model $X_{t_i} = X_{t_{i-1}} + \sigma_{t_i} U_i$ where the $t_i$ are the arrival times of a point process $N_t$. We assume that there exists a filtered probability space $(\Omega, \cF, (\cF_t)_{t\geq 0}, \bP)$, where $\cF_0$ includes all null sets and the filtration $(\cF_t)_{t\geq 0}$ is right-continuous, such that
\begin{itemize}
\item[i)] $N_t$ is a point process admitting an $\cF_t\,$-intensity $\lambda_t$ (as in  Definition D7 of Br\'{e}maud 1981); in particular $\lambda_t$ is an $\cF_t\,$-progressive process and $N_t$ is adapted to $\cF_t$;
\item[ii)] $\sigma^2_t$ is a non-negative $\cF_t\,$-predictable process; in particular $\sigma^2_t $ is $\cF_{t-}\,$-measurable;
\item[iii)] $U_i$ is $\cF_{t_i}\,$-measurable, for each $i$, with
\begin{equation*} \label{}
\bE\brac{ U_i \con \cF_{t_i-}}=0  \quad \text{and}  \quad \bE\brac{ U^2_i \con \cF_{t_i-}}=1.
\end{equation*}
\end{itemize}
\end{assu}
\noindent Note that no other condition on the distribution of $U_i$ is required. Examples for processes which fulfill these assumptions are given at the end of this chapter. The smallest filtration which satisfies the above conditions is
$$
\cF_t = \sigma\bra{  \set{N_s : s\leq t}, \set{\lambda_s : s\leq t}, \set{\sigma_s : s\leq t}, \set{U_{N_s} : s\leq t} }.
$$

\begin{prop} \label{prop_2_2}
Suppose Assumption~\ref{assu_2_1} holds. If $\sigma_t$ and $\lambda_t$ are continuous processes we have
\begin{equation} \label{ProductFormula}
\sigma_{clock}^{2}(t) = \sigma_t^{2} \cdot \lambda_t.
\end{equation}
\end{prop}

\noindent The proof can be found in Appendix~A.2.

\medskip

Thus, in the transaction-time model, the volatility can be decomposed into the product of two curves which both can be identified from the data. There exists an intuitive interpretation of this decomposition in that the formula reflects the change of time unit. Note that $\sigma^{2}_t$ can be seen as the volatility per transaction and $\sigma_{clock}^{2}(t)$ as the volatility per calendar time unit. The formula then says that
``volatility per time unit is equal to volatility per transaction multiplied by the average number of transactions per time unit''.
\begin{figure}[!htbp]
\begin{minipage}[t]{0.49\textwidth}
\includegraphics[width=\textwidth,keepaspectratio]{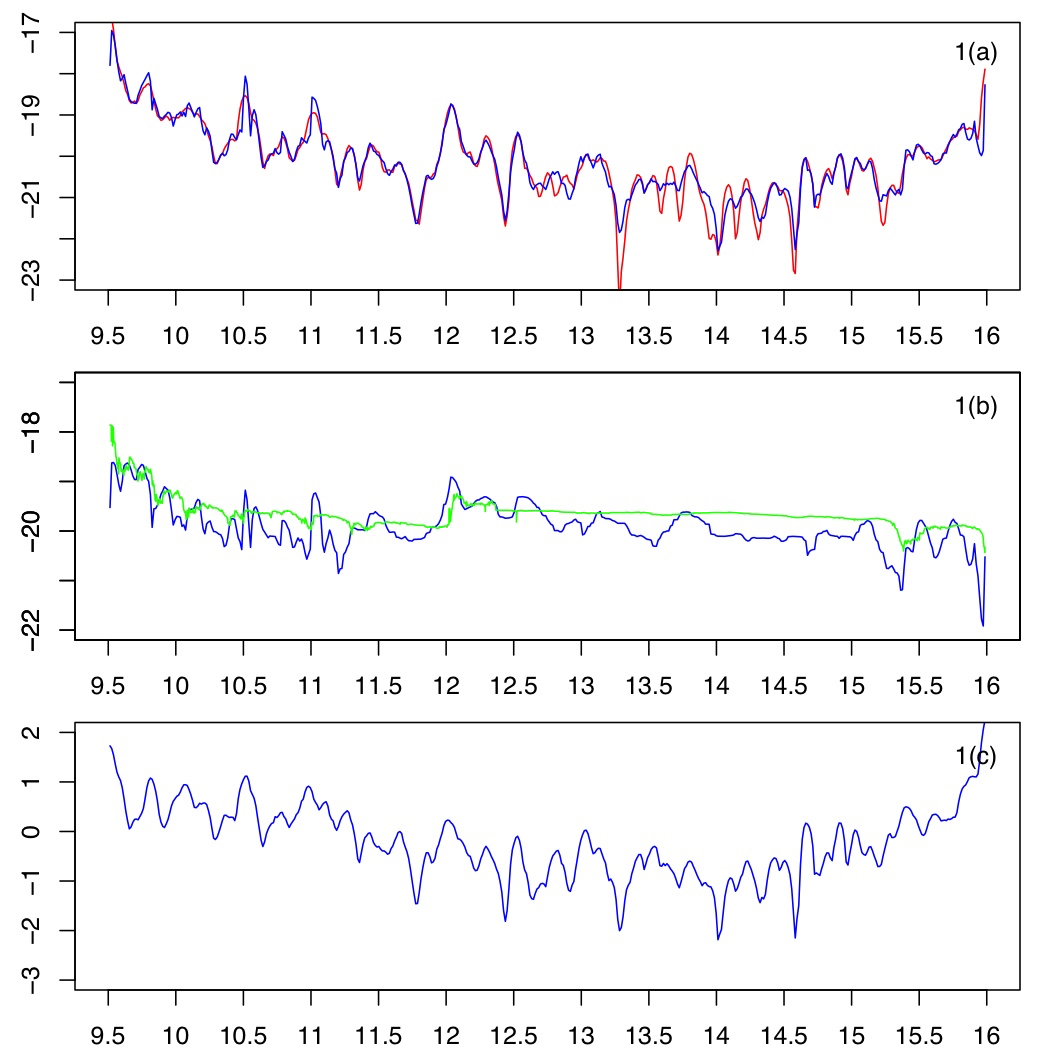}
\caption{\scriptsize Volatility of MSFT on April 1, 2014 based on 25,198 transactions, M=200: clock-time volatility(a), transaction-time volatility(b), trading intensity(c).}
\label{fig:FigureVola1b}
\end{minipage} 
\begin{minipage}[t]{0.02\textwidth}$\quad$\end{minipage}
\begin{minipage}[t]{0.49\textwidth}
\includegraphics[width=\textwidth,keepaspectratio]{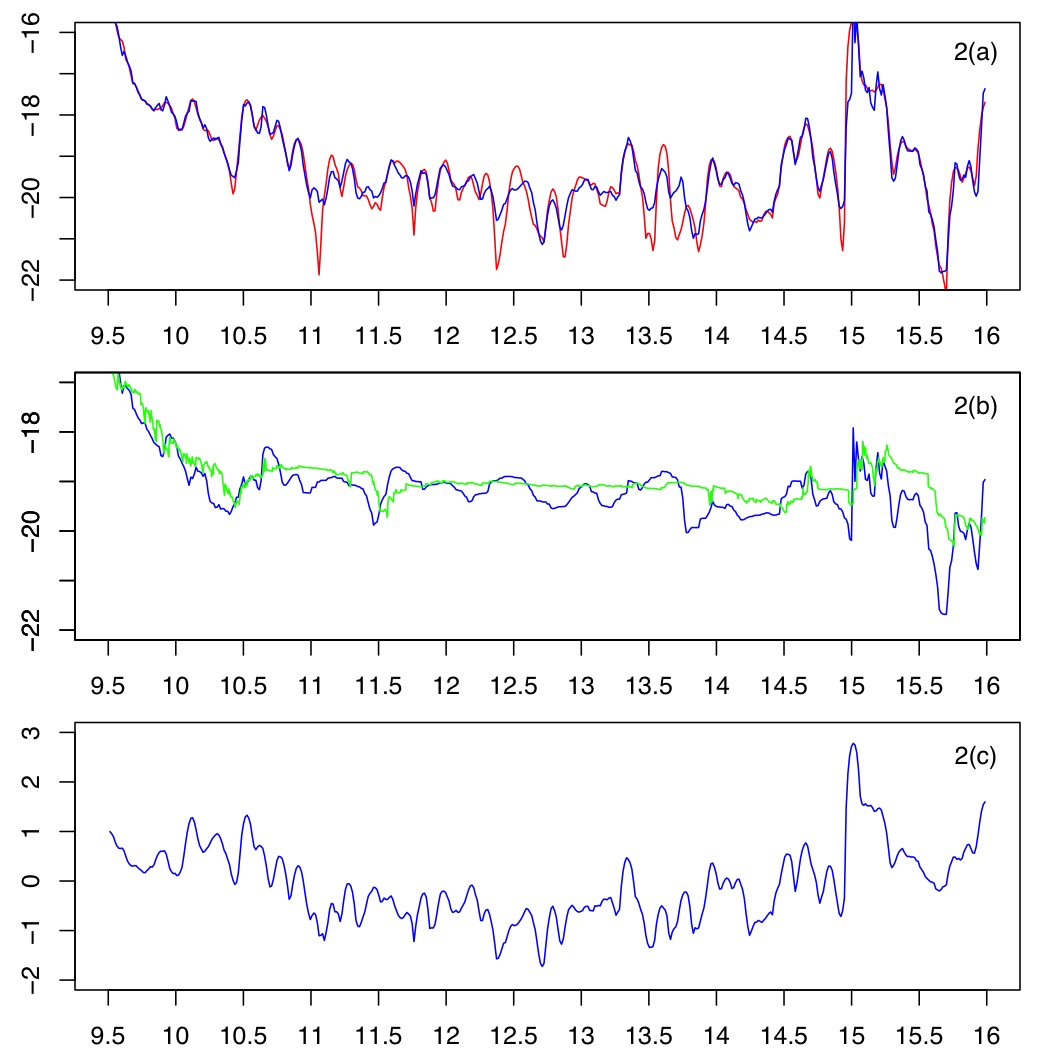}
\caption{\scriptsize Volatility of GM on April 1, 2014 based on 31,044 transactions, M=200: clock-time volatility(a), transaction-time volatility(b), trading intensity(c).}
\label{fig:FigureVola2b}
\end{minipage}
\begin{minipage}{0.49\textwidth}
\vspace*{0.3cm}
\includegraphics[width=\textwidth,keepaspectratio]{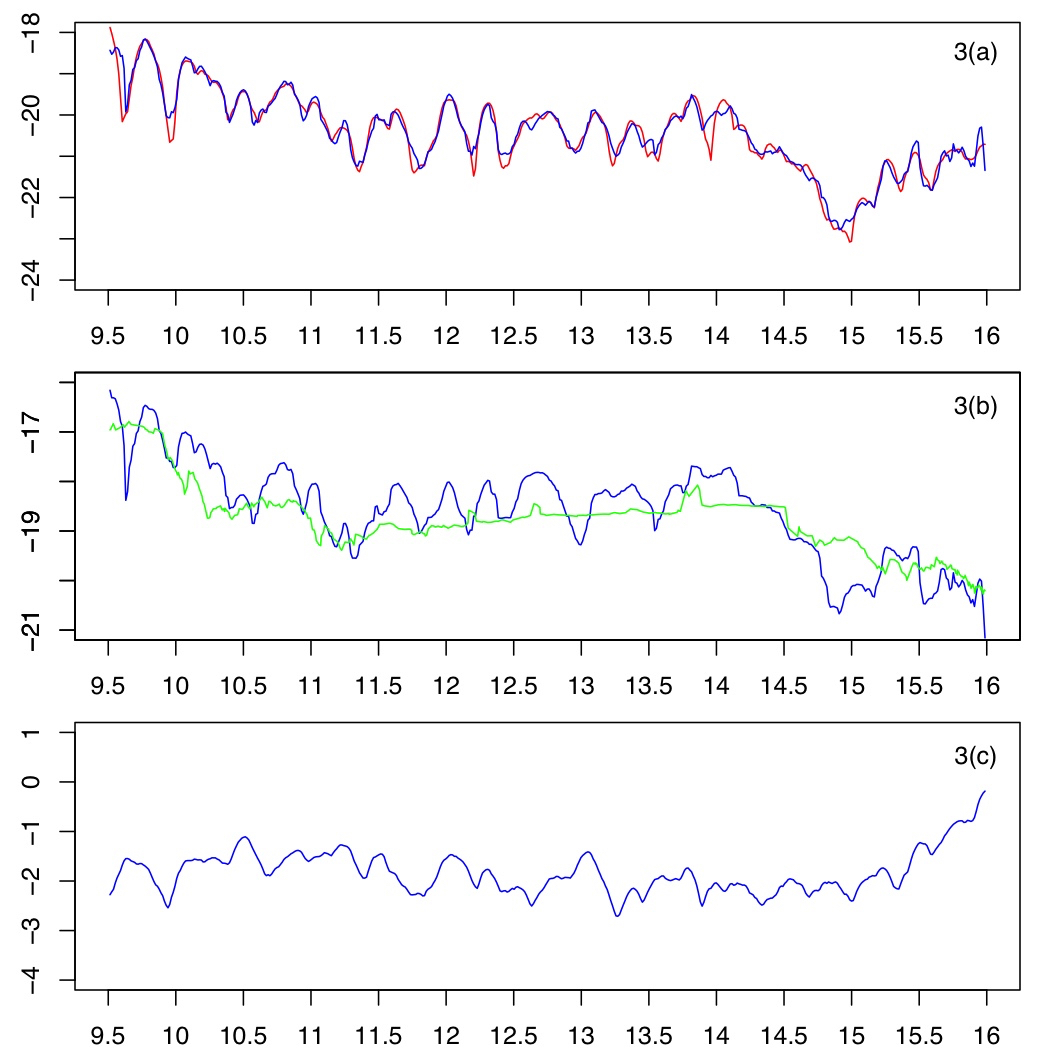}
\caption{\scriptsize Volatility of HON on April 1, 2014 based on 4,162 transactions, M=300: clock-time volatility(a), transaction-time volatility(b), trading intensity(c).}
\label{fig:FigureVola3b}
\vspace*{0.3cm}
\end{minipage} 
\begin{minipage}[t]{0.02\textwidth}$\quad$\end{minipage} 
\begin{minipage}{0.49\textwidth}
\vspace*{0.3cm}
\includegraphics[width=\textwidth,keepaspectratio]{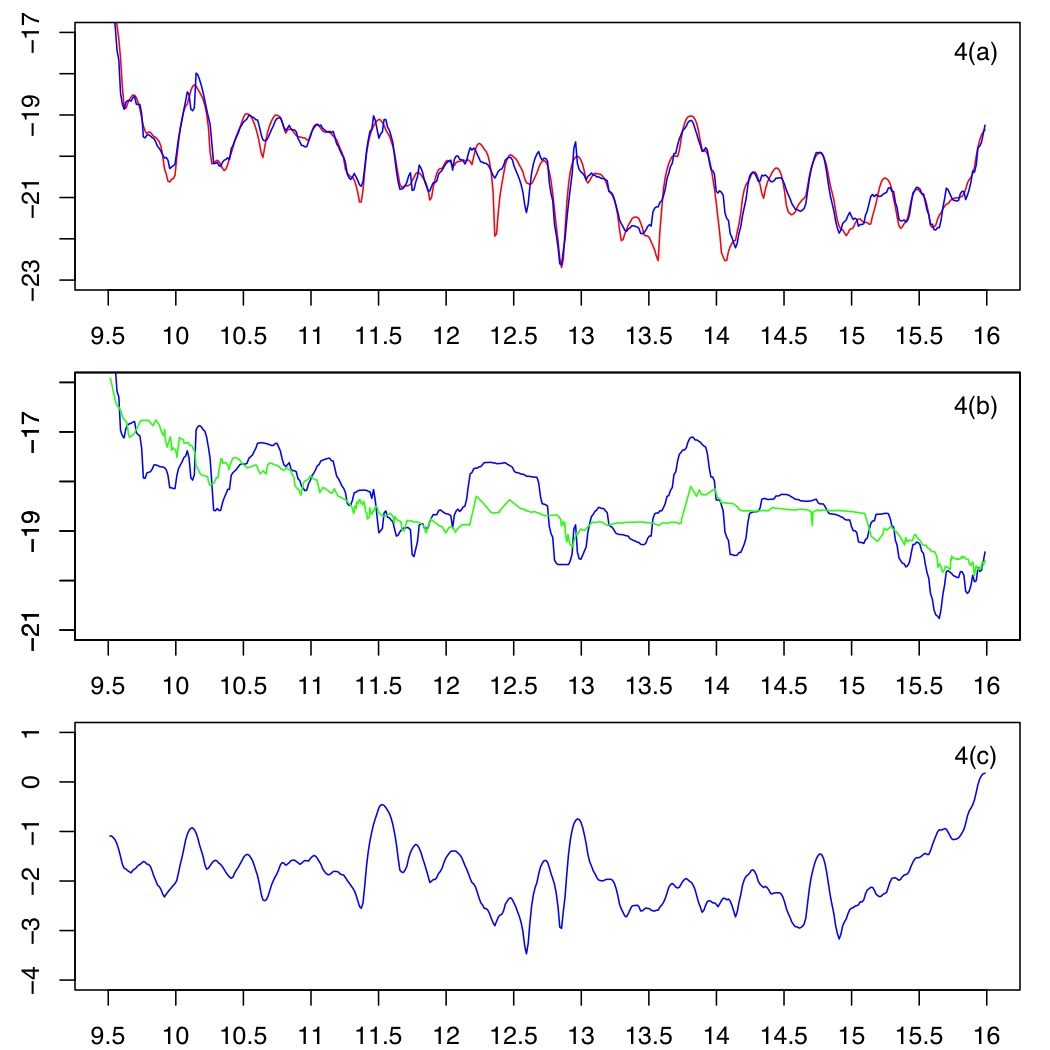}
\caption{\scriptsize Volatility of NKE on April 1, 2014 based on 4,341 transactions, M=300: clock-time volatility(a), transaction-time volatility(b), trading intensity(c).}
\label{fig:FigureVola4b}
\vspace*{0.3cm}
\end{minipage}
{Figure$\,$1--Figure$\,$4: \scriptsize The first row in each plot shows the log of clock-time volatility $\,\log \widehat{\sigma}^2_{clock,pavg}(t)$ (red) and $\,\log \widetilde{\sigma}^2_{clock,pavg}(t)$ (blue), the second row in each plot shows the log of tick-time volatility $\,\log \widehat{\sigma}^2_{pavg}(t)$ (blue) and $\,\log \widehat{\widehat{\sigma}}^2 (t)$ (green), and the third row the log trading intensity $\,\log \widehat{\lambda}(t)$ (blue). The blue estimators from row~(b) and (c) sum up to the blue estimator in (a).}
\end{figure}

Both curves $\sigma_t^{2}$ and $\lambda_t$ can be estimated by various estimates  $\widehat{\sigma}^2(t)$ and $\widehat{\lambda}(t)$. We use these estimates in two ways:

(i) to construct an alternative estimator of $\sigma_{clock}^{2}(t)$ via
\begin{equation*} \label{ProductFormulaEstimate}
\tilde{\sigma}^2_{clock}(t) \ := \ \widehat{\sigma}^2(t)\cdot\widehat{\lambda}(t);
\end{equation*}

(ii) to look at the two curves individually in order to gain more insight about the cause\linebreak \hspace*{0.5cm} and the structure of volatility. \\ \hspace*{1.1cm}

\noindent In this paper we use a kernel estimate for $\lambda_t$ (with $\int_\bR \fK (x) dx =1$ and $\fK (x)=0$ for $|x| \geq 1$ - the same for the kernels $k$ and $K$ from below)
\begin{equation} \label{IntensityEst}
\widehat{\lambda}(t_0) := \frac{1}{M} \sum^{N_T}_{i=1}  \fK\!\bra{ \frac{t_i - t_0}{M}}
\end{equation}
where $t_0 \in (0,T)$, and in order to handle microstructure noise, the pre-averaging technique of Podolskij and Vetter~(2009), extended by Jacod et al.~(2009), for the estimation of $\sigma_t^{2}$ (adapted to the present model)
\begin{align*}
\widehat{\sigma}^2_{pavg}(t_0) :=& \frac{1}{mH}\frac{1}{g_2} \sum^{i_0+m}_{i=i_0-m} k\bra{  \!\frac{i-i_0}{m}\!} \bra{\ot{Y}_{t_i}}^2 \\
& \quad -  \frac{1}{2mH} \frac{\sum^{H-1}_{l=1} h^2(\!\frac {l} {H}\!)}{g_2} \sum^{i_0+m}_{i=i_0-m} k\bra{\!\frac{i-i_0}{m}\!} \bra{Y_{t_i} - Y_{t_{i-1}}}^2
\end{align*}
where $i_0 := \inf \set{i : t_i \geq t_0}$,
\begin{align*}
\overline{\triangle Y}_{t_i} :=& \sum^{H-1}_{l=1}  g\!\bra{\frac{l}{H}} \bra{Y_{t_{i+l}} - Y_{t_{i+l-1}}}
\end{align*}
and $h(l/H) := g((l+1)/H) - g(l/H)$. $H$ is the smoothing parameter in the pre-averaging step (for more details see Section~\ref{ch:InfillAsymptotics}). This leads to the new alternative clock-time estimator based on the volatility decomposition (\ref{ProductFormula})
\begin{equation} \label{Alternative_Clock_Volatility}
\tilde{\sigma}^2_{clock,pavg}(t_0) := \ \widehat{\sigma}^2_{pavg}(t_0)\cdot\widehat{\lambda}(t_0)
\end{equation}

\noindent while the ``classical'' pre-averaging clock-time volatility estimator is
\begin{align*}
\widehat{\sigma}^2_{clock,pavg}(t_0) =&  \frac{1}{MH} \frac{1}{g_2}\sum^{N_T}_{i=1} K\!\bra{ \frac{t_i - t_0}{M}} \bra{ \overline{\triangle Y}_{t_i} }^2  \notag\\
& \quad - \frac{1}{2MH} \frac{\sum^{H-1}_{l=1} h^2(l/H)}{ g_2} \sum^{N_T}_{i=1} K\!\bra{ \frac{t_i - t_0}{M}} \bra{Y_{t_i} - Y_{t_{i-1}}}^2.
\end{align*}

In Figures~\ref{fig:FigureVola1b}--\ref{fig:FigureVola4b} we have applied both estimates to high-frequency data, approximately 4,000--30,000 transactions per day, from the NASDAQ stock exchange (MSFT = Microsoft, GM = General Motors, HON = Honeywell, NKE = Nike) - more details about the data can be found in Appendix~A.1. The first row always shows the logarithm $\,\log \widehat{\sigma}^2_{clock,pavg}(t)$ (red) and $\,\log \tilde{\sigma}^2_{clock,pavg}(t)$ (blue), the second row  $\,\log \widehat{\sigma}^2_{pavg}(t)$ (blue) and another tick-time volatility estimator $\,\log \widehat {\widehat{\sigma}}^2 (t)$ (green) from Dahlhaus and Neddermeyer (2013) - see the discussion below, and the third row the log trading intensity $\,\log \widehat{\lambda}(t)$ (blue). More details about the estimators can also be found in Appendix~A.1. Due to the additive relation
\begin{equation*} \label{}
\log \tilde{\sigma}^2_{clock,pavg}(t) = \ \log \widehat{\sigma}^2_{pavg}(t) + \log \widehat{\lambda}(t),
\end{equation*}
the blue curves in the second and third row sums up to the blue estimator in the first row. From the figures it can be seen that

\medskip

\noindent (i) row~(a) shows that the new estimator of this paper based on relation (\ref{ProductFormula}) (blue curve) nicely coincides with the classical clock-time estimator (red). This blue estimator is the sum of the blue estimators in row~(b) and (c) (in log scale);\\[8pt]
(ii) the tick-time volatility estimator $\log \widehat{\sigma}^2_{pavg}(t)$ in row~(b) is in general smoother than the clock-time estimator in row~(a) and $\log \widehat{\lambda}(t)$ in row~(c) - i.e.~the fluctuation of trading intensity in row~(c) is the major source of fluctuation of clock-time volatility in row~(a). This effect is quite clear for the high-liquid stocks MSFT and GM in Figures~\ref{fig:FigureVola1b}--\ref{fig:FigureVola2b} respectively (where $M=200$; $m=215$ and $m=265$ respectively - see Section~A.1). For the less liquid stocks HON and NKE in Figures~\ref{fig:FigureVola3b}--\ref{fig:FigureVola4b} (where $M=300$; $m=53$ and  $m=55$) the effect is less visible - in our opinion due to the considerably lower tick time bandwidth $m$.\\[8pt]
\noindent (iii) the decomposition allows to a certain extent to determine the source of volatility changes: for example the peak in 4(a) at time 13.7 is due to a peak of tick-time volatility (i.e.~most likely due to some company related news) while the peak in 2(a) at time 15.1 is mainly  due to a peak of trading intensity (i.e.~most likely due to some general - not company related - news). Similarly the decrease of volatility in 3(a) after 14.1 is company related;\\[8pt]
\noindent (iv) in particular the curves in 1(a) and 2(a) exhibit the typical U-shape over the trading day. It is notable that this U-shape is mainly a feature of the trading intensity in row~(c). In Fig.~3 and  4 the U-shape in (c) is compensated by a decrease of tick-time volatility in (b) at the end of the trading day.

\medskip

In most of our examples tick-time volatility in row~(b) is considerably smoother than trading intensity in row~(c). Beyond interpretation this has also an important consequence for estimation: microstructure noise only affects the smoother curve (b) and not (c), i.e.~coping with microstructure noise becomes easier since we may choose a larger bandwidth with effectively more data than with the classical estimator (red curve in (a)). Mathematically this is reflected in a higher rate of convergence of the estimate (see Section~\ref{ch:InfillAsymptotics} and in particular Table 1) -  in particular we may even outperform the lower bound in diffusion models  (provided that the time-rescaled model of this paper is correct).

There is an open issue about the quality of the pre-averaging estimator in (b) (a detailed investigation of this problem is beyond the scope of this paper - we just mention it briefly): The green estimator in (b) permanently is even smoother than the blue pre-averaging estimator. This estimator is completely different: it uses a nonlinear microstructure noise model with particle filtering and adaptive bandwidth selection - see Dahlhaus and Neddermeyer (2013) for details. We  have studied the behavior of both estimators by a simulation in the case where (i) the true tick-time volatility is constant (Figure~5) - emulating the case where the true curve is similar to the green estimate in \ref{fig:FigureVola1b}(b), and in the case where (ii) the true tick-time volatility is oscillating (Figure~6) - emulating the case where the true curve is similar to the blue estimate in \ref{fig:FigureVola3b}(b). As the microstructure noise model we have chosen in the simulations additive noise plus rounding. The plots in Figures~5 and 6 indicate that the green estimator resembles in particular a constant curve in a better way, and is not close to constant if the true curve is not constant, meaning that the true unknown tick-time volatility-curve in Figures~\ref{fig:FigureVola1b}(b)--\ref{fig:FigureVola4b}(b) is likely to be closer to the green curve than to the blue curve.  In particular this confirms that tick-time volatility is usually smoother than clock-time volatility and trading intensity.

\begin{figure}[!htbp]
\begin{center}
\includegraphics[width=0.7\textwidth,keepaspectratio]{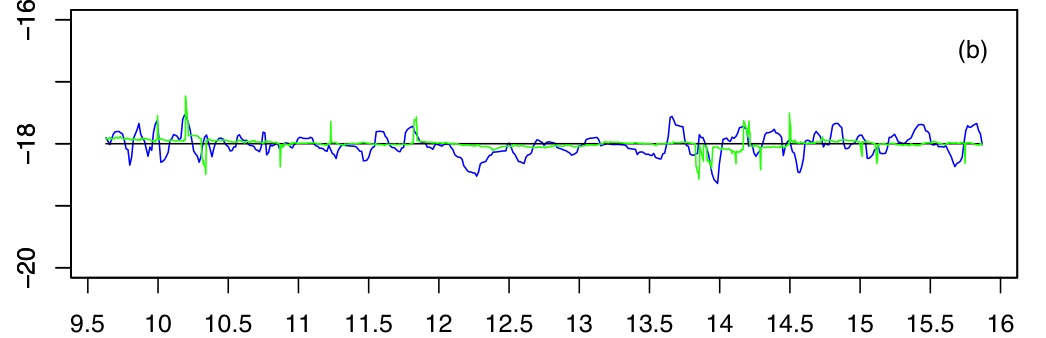}
\caption{\scriptsize Tick-time volatility estimates $\,\log \widehat{\sigma}^2_{pavg}(t)$ (blue) and $\,\log \widehat{\widehat{\sigma}}^2 (t)$ (green) for simulated transactions with true constant volatility (black)}
\end{center}
\label{fig:FigureSimulation1}
\end{figure}
\begin{figure}[!htbp]
\begin{center}
\includegraphics[width=0.7\textwidth,keepaspectratio]{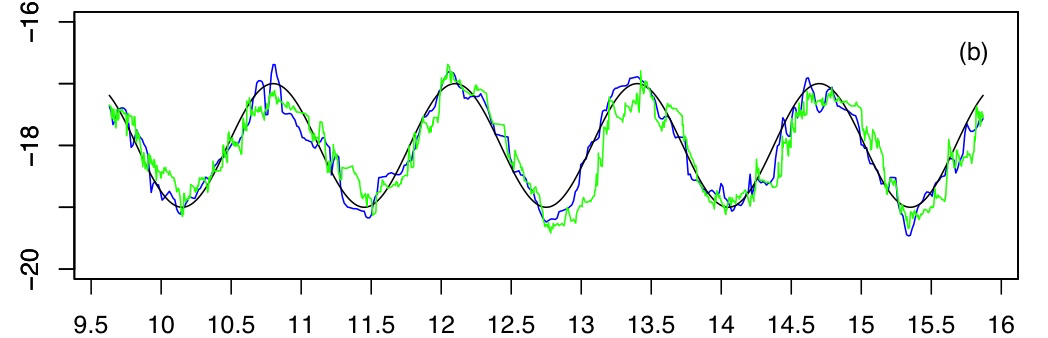}
\caption{\scriptsize Tick-time volatility estimates $\,\log \widehat{\sigma}^2_{pavg}(t)$ (blue) and $\,\log \widehat{\widehat{\sigma}}^2 (t)$ (green) for simulated transactions with true oscillating volatility (black)}
\end{center}
\label{fig:FigureSimulation2}
\end{figure}

We investigate the properties of the alternative clock-time estimator $\tilde{\sigma}^2_{clock}(t) := \widehat{\sigma}^2(t)\cdot\widehat{\lambda}(t)$ in the next section. At the end of this section we give some examples where Assumption~\ref{assu_2_1} is fulfilled.

\begin{exa} \label{Example_2_1}
\noindent (i) The simplest case is the model in \eqref{TransactionTime-Model} where $\sigma_t$ and $\lambda_t$ are deterministic and $W(\cdot)$ and $N_{\cdot}$ are independent. Even in this case the derivation of asymptotic results for the estimates is non-standard since the classical asymptotic setting cannot be applied. We therefore introduce a type of infill asymptotics for this setup in the next section.\\[8pt]
\noindent (ii) The more general model allowing for stochastic parameters is $dX_t = \sigma_t dL_t$, where $L_t$ is a pure jump process of the form $\sum^{N_t}_{i=1} U_i$. The point process $N_t$ and the sequence of innovations $U_i$ need to satisfy the conditions given in Assumption~\ref{assu_2_1}, i.e.~$N_t$ has an $\cF_t$-intensity $\lambda_t$ and $U_i$ has a martingale difference structure. Thus, $L_t$ can be seen as a generalization of a L\'evy process without a diffusion part. The leverage effect between all processes can be constructed, for example, by setting $\lambda_t := \alpha^2_t$ with
$$
d\alpha_t = \ a_tdB_t + a'_tdB'_t \quad \text{and} \quad d\sigma_t = \ b_t dB_t +  b''_t dB''_t,
$$
where $B_t$, $B'_t$, and $B''_t$ are three different $\cF_t$-Brownian motions; $a_t$, $a'_t$, $b_t$ and $b''_t$ are $\cF_t$- adapted processes. The dependence between the price process and the intensity (also the tick-time volatility) lies in the process $a_t$ (also in $a'_t$, $b_t$ and $b'_t$). For instance, we could model $a_t$ by $da_{t+} = a^*_t \ dX_t$ where $a^*_t$ is an adapted predictable process.\\[8pt]
\noindent (iii) It is interesting that GARCH-models ``almost'' fit into this framework: For example let $N_t$ be a doubly stochastic Poisson process with right-continuous $\lambda_t$ and
%
\begin{equation*} \label{}
\sigma^2_s : = a_0 +  a \set{ X_{t_{j(s)} } -X_{t_{j(s)-1} }}^2 +  b \sigma^2_{t_{j(s)}} \quad \in \cF_{s-},
\end{equation*}
where $j(s) := \max \set{ j  \ : \ t_j < s }$ and $a_0, a, b$ are positive constants. Clearly, this volatility function $\sigma^2_t$ is a left-continuous $\cF_t$-adapted step function. Since it is not right continuous we obtain instead of (\ref{ProductFormula}) with the same proof
\begin{equation*} \label{}
\sigma_{clock}^{2}(t+) = \sigma_{t+}^{2} \cdot \lambda_t
\end{equation*}
(note that $(\cF_t)_{t\geq 0}$ is right-continuous).

\medskip

\noindent (iv) In the spirit of the last example we may construct more complex examples where $\lambda_t$ is also a left-continuous $\cF_t$-adapted depending both on past arrival times of $N(t)$ (e.g.~via a similar structure as in Hawkes-models) and on past log-prices. $\sigma_t$ may in addition to the GARCH-structure from above also depend on the intensity of the point process. In that way we may explicitly model the dependence between the price process and the trading intensity. 
\end{exa}

\section{INFILL ASYMPTOTICS}\label{ch:InfillAsymptotics}

For asymptotic investigations of the estimators (e.g.~of $\tilde{\sigma}^2_{clock,pavg}(t_0) $ in \eqref{Alternative_Clock_Volatility}) we need a setup where the number of trading times increases in the neighborhood of $t_0$. A first idea would be to use a point process with intensity $n \lambda(t)$ where $n \rightarrow \infty$ (cf.~Bibinger et al. 2016, eq.~(1.2)). Instead we use here an infill asymptotic approach similar to nonparametric regression or locally stationary time series (Dahlhaus 1997). The reason for doing so is explained in Remark~\ref{rem_3_1} below.

We assume that the asset log-price process is of the form
\begin{equation} \label{RescaledTime-Model}
dX_{t,T} = \sigma\bra{\frac{t}{T}} \frac{1}{\sqrt{T}} \, dW_{N_{t,T}}, \quad \text{ for } t \in [0,T].
\end{equation}
$W_\cdot$ is a standard Brownian motion and $N_{t,T}$ is a nonhomogeneous Poisson process (NHPP) with a continuous non-negative real-valued intensity function $\lambda(t/T)$; $W_{N_{t,T}}$ is therefore a time-changed Brownian motion. In this model $\sigma^2(\cdot)$ is a real-valued deterministic continuous function called \textit{tick-time volatility}, since it responds to price variation from one trade to the next. By $t_{i,T}:= \inf\set{t : N_{t,T} \geq i}$ we denote the arrival times of the point process. In order to not complicate the notation we often avoid the double subscript, i.e.~we set $t_{i}:=t_{i,T}$ and $N_T := N_{T,T}$.
We include microstructure noise into our model (cf.~Zhang et~al. 2005; Bandi and Russell~2008); for simplicity we assume that the observations follow the linear model
\begin{equation} \label{AdditiveNoiseModel}
Y_{t_i,T} = X_{t_i,T} + \ep_i \quad \text{ for } i=1,...,N_T,
\end{equation}
where $\ep_i$ is assumed to be i.i.d. and independent of $X_{t,T}$ with
$$
\bE[\ep_i] =0, \ \var[\ep_i] = \omega^2 < \infty \ \text{and } \ \var[\ep^2_i] = \theta\omega^4, \text{ for } \theta \in \bR^+.
$$
In fact, it has been empirically shown that the assumption about the independence of the noise is justifiable for high-frequency intraday data with transaction sampling, but not for other sampling schemes such as quotation sampling, business-time sampling, 1-minute sampling, etc. (Hansen and Lunde~2006; Griffin and Oomen~2008). Throughout this section we will work under the following conditions.
\begin{assu} \label{assu_3_1}
\begin{itemize}
\item[i)] The processes $N_{\cdot,T}$ and $W_{\cdot}$ are independent;
\item[ii)] $\sigma^2(\cdot)$ is deterministic, lies in the H\"older class $\cC^{m,\gamma}[0,1]$, for $m=0,1,2$ and $0<\gamma<1$, and is bounded away from zero uniformly in $u$;
\item[iii)] $\lambda(\cdot)$ is deterministic, lies also in the H\"older class $\cC^{m',\gamma'}[0,1]$, for $m'=0,1,2$ and $0<\gamma'<1$, and is bounded away from zero uniformly in $u$.
\end{itemize}
\end{assu}
\noindent To recall the definition of the H\"older class, $f \in\cC^{m,\gamma}[0,1]$ for $0<\gamma < 1$ and $m \geq 0$, if
$$
|f^{(m)}(x + \delta)-f^{(m)}(x) | \ \leq C \cdot |\delta|^\gamma,  \ \text{for } |\delta|\to 0 \text{ and a constant } C.
$$

For model \eqref{RescaledTime-Model}, the variance of the price increment over $[t_o, t_o+bT]$ is given by
\begin{equation*}
\bE \brac{X_{t_o+bT,T} - X_{t_o,T}}^2 = \int^{u_o+b}_{u_o} \sigma^2(u)\lambda(u)\,du,
\end{equation*}
where $u_o:=t_o/T$. This implies
\begin{equation}\label{RescaledVolatilityDecom}
\sigma^2_{clock}(u_o) := \ \lim_{b \rightarrow 0} \frac{1}{b} \bE \brac{X_{t_o+bT,T} - X_{t_o,T}}^2  =  \sigma^2(u_o)\cdot\lambda(u_o),
\end{equation}
i.e.~we obtain the same volatility decomposition as in Proposition~\ref{prop_2_2} also in the rescaled model.

\medskip

Many common parametric volatility models in the literature, such as the Heston, GARCH, CIR, etc., are stochastic processes driven by a Brownian motion. Hence their realizations, being non-differentiable, lie in this class with $m=0$ and $\gamma < 1/2$. In our setting the volatility curve $\sigma^2_{clock}(\cdot)$ acts like the one that is driven by a Brownian motion whenever either of the components of \eqref{RescaledVolatilityDecom} lies in the H\"older class $\cC^{0,\gamma}$ with $\gamma < 1/2$, while the other component is allowed to be smoother (as a consequence, one may choose larger bandwidth for estimating this curve - see below).

Recently other parametric volatility models relying on a fractional Brownian motion have also been applied in price modeling. In fact, almost-all trajectories of this process with Hurst index $H \in (0,1)$ lie in the H\"older class with $m=0$ and $\gamma < H$. From this viewpoint the H\"older smoothness class seems to be an appropriate choice of smoothness class for our parameter curves.

\begin{rem} \label{rem_3_1}
\begin{itemize}
\item[i)] As mentioned above, an alternative model to (\ref{RescaledTime-Model}) seems to be a non-rescaled price model with intensity $n \lambda(t)$ where $n \rightarrow \infty$. Both approaches however are similar: In this alternative model a segment about $t_o$ of length $b$ (bandwidth of kernel estimators) would contain in the average $b \times n \lambda(t_o)$ data points while in the rescaled model of this paper a segment about $u_o:=t_o/T$ of length $b$ would contain in the average $b T \times \lambda(u_o)$ data points. Therefore, we conjecture that the two asymptotic approaches are equivalent for $n=T$. The reason for choosing the rescaled infill asymptotic model is that this approach seems to be more flexible towards future generalizations where also dependence between the price process and the trading intensity is included (which seems to be difficult if the trading intensity converges to infinity).
\item[ii)] Another aspect which looks strange at first sight is the factor $1/\sqrt{T}$ in the price model (\ref{RescaledTime-Model}) meaning that the `true tick-time volatility' is not $\sigma^{2}(t/T)$ but $\sigma^{2}(t/T)/T$. Similarly the `true clock-time volatility is $\sigma^2_{clock}(u_o)/T$, i.e.~the decomposition formula still holds regardless of this rescaling. The main reason for introducing this factor is that then the influence of microstructure noise is asymptotically the same as in the classical case:\\[6pt]
    To be more precise, consider for example the simple realized volatility estimate $\hat{v} := \frac {1} {b}\sum_{t_i \in I}  [Y_{t_i} - Y_{t_{i-1}}]^{2}$ on a segment $I$ of length $b$ in a model with constant $\sigma (\cdot) \equiv \sigma$, and trading times $t_i$ from an independent Poisson process with time-constant $\lambda (\cdot) \equiv \lambda$. We then have:\\[6pt]
    \textbf{-}\; in the classical diffusion model $dX_{t} = \sigma \, dW_t$ with $I:=[t_o, t_o+b]$,  trading times $t_i=i/T$,\\
    \hspace*{0.2cm} or alternatively trading times $t_i$  generated from a Poisson process with intensity $\lambda=T$:
    \begin{equation*} \label{}
     \bE \brac{\hat{v}} = \ \bE \frac {1} {b} \sum_{t_i \in I} \big[\sigma^{2} (t_{i} - t_{i-1}) + 2 \omega^{2}\big] \approx \ \sigma^{2} + 2 T \omega^{2}\,;
    \end{equation*}
    \textbf{-}\; in the tick-time model (\ref{RescaledTime-Model}) with factor $1/\sqrt{T}$ and $I:=[t_o/T, t_o/T+b]$, $\lambda=1$,
    \begin{align*} \label{}
    \bE \brac{\hat{v}} = \ \bE \frac {1} {b} \sum_{t_i/T \in I} \big[\sigma^{2} / T + 2 \omega^{2}\big] &\approx \frac {1} {b}\sigma^{2} b \lambda + 2 \lambda T \,\omega^{2} = \ \sigma^{2}  + 2 T \omega^{2}\,;
    \end{align*}
    \textbf{-}\; in the tick-time model (\ref{RescaledTime-Model}) without factor $1/\sqrt{T}$ and $I:=[t_o/T, t_o/T+b]$, $\lambda=1$,
    \begin{align*} \label{}
    \bE \brac{\hat{v}} = \ \bE \frac {1} {b} \sum_{t_i/T \in I} \big[\sigma^{2} + 2 \omega^{2}\big] &\approx \ \sigma^{2} \lambda T + 2 \lambda T \,\omega^{2} = \ \sigma^{2}  T + 2 T \,\omega^{2}\,,
     \end{align*}
    \hspace*{0.2cm} meaning that we had to replace $\hat{v}$ by $\hat{v}/T$ leading to $\bE \hat{v} \approx \sigma^{2} \lambda + 2 \lambda \,\omega^{2}$, i.e.~the~micro\-structure noise term were asymptotically of a lower order than in the classical case.
    
    Thus we need the factor $1/\sqrt{T}$ in model (\ref{RescaledTime-Model}) in order to make the microstructure noise problem comparable to the classical case.
\item[iii)] The independence between $N_{\cdot,T}$ and $W_{\cdot}$ enables us to make arguments conditional on the entire process $N_{\cdot,T}$ which greatly simplifies the proofs. In particular, for consecutive arrival times $t_{j-1}$ and $ t_j$ we have
\begin{equation*}
X_{t_j,T}-X_{t_{j-1},T} \stackrel{\text{law}}{=} \ \sigma\bra{ \frac{t_{j}}{T}} \frac{1}{\sqrt{T}} U_j,
\end{equation*}
since $W_{N_{t,T}}$ has the same law as $\sum^{N_{t,T}}_{i=1} U_i$, where $U_i$ are i.i.d. normally distributed variables with zero mean and unit variance and independent of $N_{\cdot,T}$.

A leverage effect between $W_\cdot$ and $N_{\cdot,T}$ could be allowed for in order to explain the correlation between market activities and price processes; see the formation of general models in Assumptions~\ref{assu_2_1}. A thoroughly asymptotic investigation of this complex transaction-time model under those assumptions is beyond the scope of this paper.
\end{itemize}
\end{rem}

We now investigate the properties of the estimates of Section~\ref{ch:voladecomposition}. For the sake of clarity we redefine them in the new infill asymptotic approach. The transaction intensity $\lambda(\cdot)$ is estimated by
\begin{equation} \label{IntensityEstimator}
\widehat{\lambda}(u_o) := \frac{1}{\fb T} \sum^{N_T}_{i=1}  \fK\bra{ \frac{t_i - u_oT}{\fb T}}  = \frac{1}{\fb T} \int^T_0 \fK\bra{ \frac{t - u_oT}{\fb T}} dN_{t,T},
\end{equation}
where the bandwidth $\fb=\fb (T) \rightarrow 0$ and $\fb T \rightarrow \infty$ as  $T \rightarrow \infty$. The kernel function $\fK$ (and also $K$ and $k$ given below) satisfies the following conditions:

\medskip

\noindent \textbf{Condition (K)} The kernel function $\fK :\bR \to\bR^+$ is a continuous, symmetric function such that $\fK (x)=0$ for $|x| \geq 1$ and $\int_\bR \fK (x) dx =1$.

\begin{thm} \label{Thm:IntensityEstimator}
Let Assumption~\ref{assu_3_1}  holds, then $\widehat{\lambda}(u_o) \xrightarrow{\bP} \lambda(u_o)$ and
$$
\sqrt{\fb T}\bra{ \widehat{\lambda}(u_o) - \bE\widehat{\lambda}(u_o) } \xrightarrow{\cD} \cN\bra{ 0, \lambda(u_o) \int_{\bR}\fK^2(x)dx}
$$
as $T \rightarrow \infty$, for $u_o \in (0,1)$. Moreover, if $\fb^{2(m'+\gamma')+1}T = o(1)$, then
$$
\sqrt{\fb T}\bra{ \bE\widehat{\lambda}(u_o)  -\lambda(u_o) - \frac{\fb^2}{2} \lambda^{(2)}(u_o)\int_{\bR} x^2\fK (x)dx \cdot I_{\set{m' =2}} } = o_p(1),
$$
particularly
\begin{align*}
&\sqrt{\fb T}\bra{ \widehat{\lambda}(u_o) -\lambda(u_o) - \frac{\fb^2}{2} \lambda^{(2)}(u_o)\int_{\bR} x^2\fK (x)dx \cdot I_{\set{m' =2}} } \xrightarrow{\cD} \cN\bra{ 0, \lambda(u_o) \int_{\bR}\fK^2(x)dx}.
\end{align*}
\end{thm}

As discussed in Section~\ref{ch:voladecomposition} we compare a classical pre-averaging estimator for $\sigma^2_{clock}(u_o)$ with a new estimator based on the decomposition formula (\ref{RescaledVolatilityDecom}). The classical estimator in the infill-framework is defined by
\begin{align}\label{ClockTimeVolEstimator}
\widehat{\sigma}^2_{clock,pavg}(u_o) :=& \frac{1}{bH} \frac{1}{g_2}\sum^{N_T}_{i=1} K\bra{ \frac{t_i - u_oT}{bT}} \bra{ \overline{\triangle Y}_{t_i,T} }^2  \notag\\
& \ - \frac{1}{2bH} \frac{\sum^{H-1}_{l=1} h^2(l/H)}{ g_2} \sum^{N_T}_{i=1} K\bra{ \frac{t_i - u_oT}{bT}} \bra{Y_{t_i,T} - Y_{t_{i-1},T}}^2,
\end{align}
with pre-averaging steps $\overline{\triangle Y}_{t_i,T}$ given by
\begin{align*}
\overline{\triangle Y}_{t_i,T} :=& \sum^{H-1}_{l=1}  g\bra{\frac{l}{H}} \bra{Y_{t_{i+l},T} - Y_{t_{i+l-1},T}} = -\sum^{H-1}_{l=1} h\bra{\frac{l}{H}} Y_{t_{i+l},T}.
\end{align*}
We define $h(l/H) := g((l+1)/H) - g(l/H)$, where $g$ is another differentiable weight function defined on $[0,1]$ with $g(0)=g(1) =0$ and which has piecewise Lipschitz continuous derivatives $g^{(1)}$. The kernel function $K$ and the bandwidth $b$ are similar to those of the intensity estimate, i.e.~$K$ satisfies condition~(K) and $b$ depends on the time span $T$ such that $b \to 0$ and $bT \to \infty$ as $T \to \infty$. Moreover the pre-averaging block size $H=H(b)$ also depends on $T$ such that $H \to \infty$ and $H/bT \to 0$ as $T \to \infty$. We assume that the limits $\sum^{H-1}_{l=1} g(l/H)^n /H $ and $\sum^{H-1}_{l=1} g^{(1)}(l/H)^n /H $ exist and equal $g_n:=\int^1_0 g(x)^n dx$ and $g'_n:=\int^1_0 g^{(1)}(x)^n dx$ respectively, for $n\in \bN$.

The subscript \textit{pavg} stands for the name of the procedure, \textit{pre-averaging}. As its name suggests, we first calculate the average of log returns weighted by a function $g$ over each block of size $H$ and then apply a local sum of squares of these averages (similar to the filtering of realized volatility) to construct a spot volatility estimate. By doing this, the variance of the noise is reduced by a factor of $1/H$, as can be seen in our proof (see also Jacod et al.~2009). In particular, this block size $H$ will play a crucial role in this setting along with the main bandwidth size $b$. Finally, a bias term induced by the additional measurement error will be corrected by the second term of \eqref{ClockTimeVolEstimator}.

\begin{thm} \label{Thm:ClockTimeVolatilityEstimator}
Suppose that Assumption~\ref{assu_3_1} is fulfilled. Let $K$ have bounded first derivatives and let the block size $H =\delta\cdot T^{1/2}$ for $\delta \in (0,\infty)$ and $b^{2\alpha+1 }T^{1/2} =o(1)$ with $\alpha =\min\set{m+\gamma,m'+\gamma'}$. Then
\begin{align*}
&\sqrt{bT^{1/2}} \set{ \widehat{\sigma}^2_{clock,pavg}(u_o) -\sigma^2_{clock}(u_o) - BIAS} \xrightarrow{\cD} \ \cN \bra{0, \ \delta\eta^2_A+\frac{1}{\delta}\eta^2_B + \frac{1}{\delta^3}\eta^2_C}
\end{align*}
as $T \to \infty$, for $u_o\in (0,1)$, where
\begin{align*}
\eta^2_A =& \ 2\sigma^4(u_o)\lambda(u_o)\int_\bR K^2(x)dx,  \qquad \qquad \ \eta^2_B = 4\omega^2 \sigma^2(u_o)\lambda(u_o) (g'_2/g_2)\int_\bR K^2(x)dx, \\
\eta^2_C =& \ 2\omega^4 \lambda(u_o) (g'_2/g_2)^2\int_\bR K^2(x)dx, \quad BIAS = \frac{1}{2} \bra{\sigma^2(u_o)\lambda(u_o)}^{(2)} b^2 \int_\bR x^2K(x)dx \cdot I_{\set{m=m' =2}}.
\end{align*}
\end{thm}

We see that our result based on the transaction-time model is different from the classical result based on the standard diffusion model in the way that the asymptotic bias and variance rely on the transaction intensity. The consistency of this estimate has been implicitly proven, leading to many consistent estimators for functionals of spot volatility (including the integrated volatility) by using continuous-mapping theorem. For example, the unknown component in the asymptotic variance $\sigma^4(u_o)\lambda(u_o)$ can be estimated by the square of $\widehat{\sigma}^2_{clock,pavg}(u_o)$ divided by $\widehat{\lambda}(u_o)$. Alternatively one could construct another consistent estimator for the asymptotic variance by applying the pre-averaging technique (Jacod et al. 2009, eq. (3.7)).

Due to the concept of volatility decomposition, we now formulate an alternative estimator for the clock-time volatility with the following product
\begin{equation} \label{ClockTimeVolEstimator_Alternative}
\tilde{\sigma}^2_{clock,pavg}(u_o) \ := \ \widehat{\sigma}^2_{pavg}(u_o)\cdot\widehat{\lambda}(u_o).
\end{equation}
For this product, the intensity estimate \eqref{IntensityEstimator} is used and the tick-time volatility estimator is given~by
\begin{align}  \label{TickTimeVolEstimator}
\widehat{\sigma}^2_{pavg}(u_o) :=&  \frac{T}{NH}\frac{1}{g_2} \sum^{i_o+N}_{i=i_o-N} k\bra{  \frac{i-i_o}{N}} \bra{\ot{Y}_{t_i,T}}^2 \notag\\
& \ - \frac{T}{2NH} \frac{\sum^{H-1}_{l=1} h^2(l/H)}{g_2} \sum^{i_o+N}_{i=i_o-N} k\bra{\frac{i-i_o}{N}} \bra{Y_{t_i,T} - Y_{t_{i-1},T}}^2,
\end{align}
where $i_o := \inf \set{i : t_i \geq u_o T}$, i.e.~$t_{i_o}$ is the first arrival time after or at the time point of interest $t_o$. The kernel function $k$ and the segment length $N$ satisfy Condition~(K) and $N=N(T)\to \infty$ and $N/T \to 0$ as $T \to \infty$. In fact, this estimator also relies on the pre-averaging approach, and therefore the same conditions for $g$ and $h$ are taken from the preceding clock-time volatility estimator  with the block size $H=H(N)$ satisfying $H \to \infty$ and $H/N \to 0$ as $T\to \infty$. Note that we have used the same letters $g$, $h$, and $H$ for both estimators in order to not complicate the notation.

At first sight, both estimates $\widehat{\sigma}^2_{clock, pavg}(\cdot)$ and $\widehat{\sigma}^2_{pavg}(\cdot)$ look very similar, as they are based on the filtering of pre-averaging estimators. However, there is a distinction between these two: one is based on tick time and the other is based on clock time. More precisely, in \eqref{TickTimeVolEstimator}   the (exactly) $N$-nearest observed pre-averaged terms from both sides of the considered time point $u_o = t_o/T$ are taken into account so that the influence of the arrival rate is removed, while $K$ in \eqref{ClockTimeVolEstimator} uses all pre-averaged terms inside the interval $[t_o-bT, t_o+bT]$. As a matter of fact, the number of transactions/pre-averaged terms over this interval is random and depends on the trading intensity. 

\begin{thm} \label{Thm:TickTimeVolEstimator}
Under Assumption~\ref{assu_3_1}, the pre-filtering block size $H=\delta\cdot T^{1/2}$ for $\delta \in (0,\infty)$, and $k$ has bounded first derivatives, we obtain
\begin{equation*}
\sqrt{\frac{N}{T^{1/2}}}\set{\widehat{\sigma}^2_{pavg}(u_o) - \sigma^2(u_o)} \xrightarrow{\cD} \ \cN\bra{0, \ \delta\xi^2_A + \frac{1}{\delta}\xi^2_B + \frac{1}{\delta^3}\xi^2_C}
\end{equation*}
where
\begin{align*}
&\xi^2_A = \ 2\sigma^4(u_o)\int_\bR k^2(x)dx, \ \ \xi^2_B = 4\omega^2 \sigma^2(u_o) \frac{g'_2}{g_2}\int_\bR k^2(x)dx \ \text{ and } \ \xi^2_C = 2\omega^4 \bra{\frac{g'_2}{g_2}}^2\int_\bR k^2(x)dx,
\end{align*}
under the segment conditions
\small
\begin{align} \label{eq_3_14}
N^{1+2\gamma}/T^{1/2+2\gamma} \to 0 \quad &\text{ for } m=0, \quad  &N^{3+\gamma^*}/T^{5/2+\gamma^*} \to 0 \quad \text{ for } m=1 \text{ and } m'=0, \notag\\
N^{3+\gamma}/T^{5/2+\gamma} \to 0 \quad &\text{ for } m=1 \text{ and } m'=1,2,\quad  &N^{3+\gamma'}/T^{5/2+\gamma'} \to 0 \quad \text{ for } m=2 \text{ and } m'=0, \notag\\
N^4/T^{7/2}\to 0 \quad &\text{ for } m=2 \text{ and } m'=1,2,
\end{align}
\normalsize
with $\gamma^* = \min\set{\gamma,\gamma'}$ as $T \to \infty$.
\end{thm}

It is remarkable that the bias derivation for this estimate is much more complicated than that of $\widehat{\sigma}^2_{clock, pavg}(\cdot)$ in the case of higher orders of smoothness ($\min(m,m') \geq 1$). In particular, the explicit bias term $BIAS$ in Theorem~\ref{Thm:ClockTimeVolatilityEstimator} cannot be stated, even though the kernel function is symmetric and the parameter functions are twice differentiable ($m,m'=2$). This weakness leads to a reduction in the rate of convergence in many cases, particularly when $\min(m,m') \geq 1$ (more precisely, $(I_{4,2,1})$ and $(I_{4,2,2})$ in the proof of Theorem~\ref{Thm:TickTimeVolEstimator} only vanish if $N^4/(HT^3) \to 0$). The above segment conditions are given to enable us to neglect some asymptotic bias terms.

From \eqref{ClockTimeVolEstimator_Alternative}, the alternative estimator is clearly consistent, since it is the product of two consistent estimators. The limit distribution is given below, where the rate of convergence will depend on the convergence rate $\fb T$ and $N/T^{1/2}$ in Theorems~\ref{Thm:IntensityEstimator} and \ref{Thm:TickTimeVolEstimator} respectively. For example, if $N/T^{1/2} = o(\fb T)$, we have
\begin{align*}
&\sqrt{\frac{N}{T^{1/2}}} \set{ \tilde{\sigma}^2_{clock,pavg}(u_o)  - \sigma^2_{clock}(u_o)}  \\
&\quad = \ \widehat{\lambda}(u_o)\sqrt{\frac{N}{T^{1/2}}} \set{ \widehat{\sigma}^2_{pavg}(u_o)  - \sigma^2(u_o) }  + \sigma^2(u_o) \sqrt{\frac{N}{T^{1/2}}}\frac{1}{\sqrt{\fb T}}\cdot \sqrt{\fb T} \set{ \widehat{\lambda}(u_o) - \lambda(u_o) - BIAS_{\lambda} } \\
&\qquad + \sigma^2(u_o)  \sqrt{\frac{N}{T^{1/2}}} \cdot BIAS_{\lambda} \\
&\quad = \widehat{\lambda}(u_o)\sqrt{\frac{N}{T^{1/2}}} \set{ \widehat{\sigma}^2_{pavg}(u_o)  - \sigma^2(u_o) }  +  o(1)  + \sigma^2(u_o)  \sqrt{\frac{N}{T^{1/2}}}\cdot BIAS_{\lambda}
\end{align*}
Thereby, the resulting limit distribution is dominated by the limit of $\widehat{\sigma}^2_{pavg}(\cdot)$ given in the last theorem, not that of $\widehat{\lambda}(\cdot)$.

\begin{thm} \label{Thm:AlternativeClockTimeVolEstimator}
Let all assumptions be satisfied and suppose that the bandwidth $\fb$ and the segment length $N$ fulfill the conditions given in Theorem~\ref{Thm:IntensityEstimator} and \ref{Thm:TickTimeVolEstimator} respectively. For $u_o \in (0,1)$ we obtain
\begin{equation*}\label{}
\sqrt{\frac{N}{T^{1/2}}} \set{ \tilde{\sigma}^2_{clock,pavg}(u_o)  - \sigma^2_{clock}(u_o) } \xrightarrow\cD \cN(0,V^2)
\end{equation*}
with
\small
\begin{align*}
&V^2 = \lambda^2(u_o)\set{\delta\xi^2_A + \frac{1}{\delta}\xi^2_B + \frac{1}{\delta^3}\xi^2_C}  I_{ \set{ m=m'=0, \gamma' > \frac{\gamma}{2\gamma+2} \text{ \textbf{or} } m=1,m'=0, \gamma' > \frac{\gamma^*+2}{2\gamma^*+8} \text{ \textbf{or} } m=2, m'=0, \gamma'> \frac{\sqrt{65}-7}{4} \text{ \textbf{or} }  m'=1,2   }},
\end{align*}
\normalsize
and
\begin{equation*}\label{}
\sqrt{\fb T} \set{ \tilde{\sigma}^2_{clock,pavg}(u_o)  - \sigma^2_{clock}(u_o)} \xrightarrow\cD \cN(0,W^2)
\end{equation*}
with
\small
\begin{align*}
&W^2= \ \sigma^4(u_o)\lambda(u_o)\int_\bR \fK^2(x)dx  I_{ \set{ m=m'=0, \gamma' \leq \frac{\gamma}{2\gamma+2} \text{ \textbf{or} } m=1,m'=0, \gamma' \leq \frac{\gamma^*+2}{2\gamma^*+8} \text{ \textbf{or} } m =2, m'=0, \gamma'\leq \frac{\sqrt{65}-7}{4} }} \\
&+ \ c_1 \lambda^2(u_o)\set{\delta\xi^2_A + \frac{1}{\delta}\xi^2_B + \frac{1}{\delta^3}\xi^2_C} I_{ \set{ m=m'=0, \gamma' = \frac{\gamma}{2\gamma+2} \text{ \textbf{or} } m=1,m'=0, \gamma' = \frac{\gamma^*+2}{2\gamma^*+8} \text{ \textbf{or} } m =2, m'=0, \gamma'= \frac{\sqrt{65}-7}{4}  }},
\end{align*}
\normalsize
where $c_1 := bT/(N/T^{1/2}) $ if $bT$ and $N/T^{1/2}$ are of the same order.
\end{thm}

\begin{rem}
\begin{itemize}
\item[i)] For the pre-averaging steps, the block size $H$ in \eqref{ClockTimeVolEstimator} may differ from $H$ in \eqref{TickTimeVolEstimator} (both are related to $T$). For this theoretical investigation, we select $H = O(T^{1/2})$ to balance the rate of convergence of the limit distributions (I), (II), and (III) (see the beginning of the proofs of Theorems~\ref{Thm:ClockTimeVolatilityEstimator} and \ref{Thm:TickTimeVolEstimator} in section~A.3) in order to obtain those asymptotic normality results.
\item[ii)] From a practical point of view, all of the unknown components in the asymptotic variances can be estimated by using the existing statistics presented in this section. The variance $\omega^2$ of the microstructure noise can be estimated by $\widehat{\sigma}^2_{clock}(u)/ (2T\widehat{\lambda}(u))$, where $\widehat{\sigma}^2_{clock}(u):= \sum^{N_T}_{i=1} \frac{1}{b} K\bra{ \frac{t_i - uT}{bT}} \bra{ Y_{t_i,T} - Y_{t_{i-1},T}}^2 $ is a filtered realized volatility, since $\widehat{\sigma}^2_{clock}(u) = 2T\omega^2 \lambda(u) + o_p(1)$. The data-adaptive choice  of the smoothing parameters $b$, $\fb$, $N$, $H(b)$, and $H(N)$ remains to be solved.
\end{itemize}
\end{rem}

It is well-known that the presence of microstructure noise causes a reduction in the rate of convergence of volatility estimation. The advantage of our decomposable estimator with respect to microstructure noise is, that the noise does not disturb the transaction-times but only the transaction-prices. In Section~\ref{ch:voladecomposition}, our empirical analysis suggests that the tick-time volatility curve is in general less fluctuating than the intensity curve. Therefore coping with microstructure noise becomes easier, as we may choose a larger window for $\widehat{\sigma}^2_{pavg}(\cdot)$ with effectively more data than with $\widehat{\sigma}^2_{pavg, clock}(\cdot)$, which is as rough as the intensity curve. Mathematically, this leads to a higher rate of convergence of the volatility estimator (in some cases even better than the lower bound of the estimation for spot volatility in the standard noisy model - see below).

From this point of view, we now compare the performance of the estimates $\widehat{\sigma}^2_{clock,pavg}(\cdot)$ (based on the classical pre-averaging method) and $\tilde{\sigma}^2_{clock,pavg}(\cdot)$ (based on the volatility decomposition). The same kernel functions $\fK$, $k$, and $K$ and weighting function $g$ are used in both estimates. Before discussing the comparison which is summarized in Table~\ref{table:Comparison_PAVG_Volatility}, we explicitly demonstrate one of those cases in detail (the other cases can be done similarly). Let $m=m'=0$ and $\gamma' \leq \gamma/2(\gamma+1)$. According to Theorems~\ref{Thm:AlternativeClockTimeVolEstimator} and \ref{Thm:ClockTimeVolatilityEstimator} we have
$$
\sqrt{\fb T} \set{ \tilde{\sigma}^2_{clock,pavg}(u_o)  - \sigma^2_{clock}(u_o)} \xrightarrow\cD  \cN(0,W^2)
$$
under the bandwidth condition $\fb^{2\gamma'+1}T \to 0$, and
$$
\sqrt{bT^{1/2}} \set{ \widehat{\sigma}^2_{clock,pavg}(u_o) -\sigma^2_{clock}(u_o) } \xrightarrow{\cD}  \cN(0, \delta\eta^2_A+\frac{1}{\delta}\eta^2_B + \frac{1}{\delta^3}\eta^2_C)
$$
under  $ b^{2\gamma'+1}T^{1/2} \to 0$. These constraints imply
$$
\fb T = \ o\bigl( T^{\frac{2\gamma'}{2\gamma'+1}} \bigr) \quad \text{and} \quad b T^{1/2} = \ o\bigl( T^{\frac{\gamma'}{2\gamma'+1}} \bigr),
$$
meaning that the rate of convergence of $\tilde{\sigma}^2_{clock,pavg}(\cdot)$ is much faster than that of $\widehat{\sigma}^2_{clock,pavg}(\cdot)$.

Moreover, in the standard diffusion model with the presence of noise, a lower bound for spot volatility estimation is derived in a minimax sense with respect to the $L_2$-loss function. This equals $n^{\frac{-\alpha}{2\alpha+1}}$, given a H\"older-exponent of $\alpha$ for the spot volatility function, where $n$ is the number of subdivisions (Munk and Schmidt-Hieber~2010). Our estimate $\tilde{\sigma}^2_{clock,pavg}(\cdot)$ is better than that bound in this particular case and in many other cases, see Table~\ref{table:Comparison_PAVG_Volatility}. Thus, the approach based on volatility decomposition of transaction-time models outperforms previous approaches applied to the standard model in these cases (provided that the time-rescaled model of this paper is correct).

\begin{table*}[t]
\caption {\footnotesize Comparison of rates of convergence and asymptotic variances between spot volatility estimators}
\label{table:Comparison_PAVG_Volatility}
\small
\begin{center}
\bgroup
\begin{tabular}{  c  c c  c  c  c  }
\hline \hline
    Case & \multicolumn{2}{c}{Conditions} & Rate($\tilde{\sigma}^2_{clock,pavg}$) vs. & $\var\brac{\tilde{\sigma}^2_{clock,pavg}}$ vs. \\
& & & Rate($\widehat{\sigma}^2_{clock,pavg}$) & $\var\brac{\widehat{\sigma}^2_{clock,pavg}}$ \\ \hline
  c1 & $m=0, \ m'=0$ & $\gamma > \gamma'$ &                                       			  & \\ 
  c2 & $m=1, \ m'=0$ &  - &                                      			  & \\ 
  c3 & $m=1, \ m'=1$ & $2\gamma' < \gamma$ &    faster                                  			  & \\ 
  c4 & $m=2, \ m'=0$ & - &                                       			  & \\ 
  c5 & $m=2, \ m'=1$ & $\gamma' < 1/2$ &                                       			  & \\[8pt] 

  c6 & $m=0, \ m'=0$ & $\gamma \leq \gamma'$ &                                       			  & \\ 
  c7 & $m=0, \ m'=1$ & - &      same                                 			  & smaller, if $\lambda(\cdot) <1$; \\ 
  c8 & $m=0, \ m'=2$ & - &                                     			  &  larger, otherwise.\\[8pt] 

  c9 & $m=1, \ m'=1$ & $\gamma\leq2\gamma' $ &                                         			  & \\ 
  c10 & $m=1, \ m'=2$ & - &      	              unknown                 			  &\\ 
  c11 & $m=2, \ m'=1$ & $\gamma' \geq 1/2$ &                                    			  & \\ 
  c12 & $m=2, \ m'=2$ & -  &                                       			  & \\ \hline
\end{tabular}
\egroup
\end{center}
\normalsize
\end{table*}

Table~\ref{table:Comparison_PAVG_Volatility} compares the performance of $\widehat{\sigma}^2_{clock,pavg}(\cdot)$ and $\tilde{\sigma}^2_{clock,pavg}(\cdot)$ in terms of their rates of convergence and asymptotic variances. We see that the alternative volatility estimator outperforms the standard estimator in the sense that its rate of convergence is  improved, see c1--c5. In these cases, the smoothness order of the tick-time volatility is higher than that of the clock-time volatility and intensity (this is clear from above, as we can choose larger window for estimating tick-time volatility than for estimating clock-time volatility). In fact, in the empirical examples in Section~\ref{ch:voladecomposition} the tick-time volatility curve is smoother than the transaction intensity curve. Thus, we can significantly improve volatility estimation by considering the volatility decomposition. In c6--c8, both estimators possess the same rate of convergence, but the asymptotic variance of $\tilde{\sigma}^2_{clock,pavg}(\cdot)$ turns out to be smaller if the intensity is less than 1; otherwise $\widehat{\sigma}^2_{clock,pavg}(\cdot)$ yields better results. Lastly, in cases c9--c12 the situation is unknown since it is very difficult to derive the bias of $\widehat{\sigma}^2_{pavg}(\cdot)$ explicitly in this situation when $\min(m,m')\geq 1$. The conditions needed to show that the bias is of lower order indicate however, that the alternative estimator has a slower rate of convergence than the classical one.

\section{CONCLUSION} \label{ch:conclusion}

In this paper we have advocated the use of a spot volatility estimate based on a volatility decomposition in a time-changed price-model according to the trading times. In this model clock-time volatility splits up into the product of two curves, namely tick-time volatility and trading intensity. Both curves can be identified and we have argued that both curves contain valuable information about the original volatility curve. For example U-shape and the increase of volatility at the end of the trading day are in our opinion solely features of trading intensity while the influence of company related news mainly hits tick-time volatility or both curves.

An important finding in our view is that the tick-time volatility curve is often much smoother than the clock-time volatility curve of high-liquid equities. This means that the major part of fluctuations in clock-time volatility is due to fluctuations of the trading intensity.

There is an important consequence of these findings also for statistical inference: microstructure noise does not influence the estimate of trading intensity but only the estimator of tick-time volatility. Since this usually is the smoother curve we may choose a larger bandwidth with the benefit of a faster rate of convergence and a better coping of microstructure noise. In particular, one may outperform the rate of convergence of the optimal estimator in the classical diffusion model.

For the mathematical investigation of this model we have introduced an infill asymptotic approach, and derived the asymptotic properties of the new estimator in the case of a deterministic volatility curve and a deterministic intensity curve of a point process. If both curves are replaced by stochastic processes one may use ideas along the lines of the work by Koo and Linton (2012) on locally stationary diffusion models in combination with similar models for the point process of transaction times such as in Roueff et al. (2016). A comprehensive treatment of this situation seems challenging and beyond the scope of this paper.

From an applied point of view it is also of high interest to find proper models where the trading intensity depends on past log-prices, and where the volatility has some GARCH-type structure - possibly depending in addition on the past intensity of the point process.






\section*{APPENDIX}

\subsection*{A.1  Details of the Data Analysis}  

In our empirical study we use tick-resolution trading history provided by the data vendor \textit{Quant-Quote TickView}. We analyze intraday transaction data from the NASDAQ stock exchange.

In general we clean raw data before analyzing it in the following main steps: i) deleting all pre- and after-market data, i.e.~only transactions between 09:30 AM--04:00 PM are considered; ii) filtering raw data from the outliers, such as price errors, and deleting entries with abnormal sale conditions. Although the accuracy of transaction time is down to milliseconds, the resolution of timestamps is limited to only one second, which leads to the possibility of having multiple consecutive transactions occurring at the same time. Nevertheless, the order of the trades is correctly placed. In the case of multiple trades, which is often the case for liquidly traded equities, these time points will be separated into equally-spaced times, for example $t_{10}, t_{11}$ and $ t_{12}$ occurring at time 34210 (= 09:30:10 AM) are adjusted to $t_{10}=34210$, $t_{11}=34210.33$ $t_{12}= 34210.67$. We note that there are 6.5 market hours in a trading day, which is equal to $T=23,400$ seconds.

Figures~\ref{fig:FigureVola1b}--\ref{fig:FigureVola4b} show the analysis of transaction data on April 1, 2014; MSFT = Microsoft (25,198 transactions), GM = General Motors (31,044 transactions), HON = Honeywell (4,162 transactions) and NKE = Nike (4,341 transactions) with the estimators described in Section~\ref{ch:voladecomposition} (remember that these estimators differ from the estimators from Section~\ref{ch:InfillAsymptotics} by the factor $1/T$). For the investigation of the high-liquid stocks in Figures~\ref{fig:FigureVola1b}--\ref{fig:FigureVola2b} we chose the time-bandwidth $M=200$ (the estimate uses $2M=400$), the tick-bandwidth $m= \lfloor 200\cdot \set{\text{\# of trades}} /T \rfloor$ (resulting in $m=215$ for MSFT; $m=265$ for  GM) and $H=15$, i.e. $M$ and $m$ were chosen to cover the same range (for $\widehat{\sigma}^2_{clock,pavg}(t_o)$, all observed pre-averaged terms $\overline{\triangle Y}_{t_i}$ inside the interval $(t_o-M, t_o+M]$ contribute to the estimator at time $t_o$, whereas exactly $m$ observed pre-averaged terms from the left- and the right-hand side of $t_o$ contribute to the estimator $\widehat{\sigma}^2_{pavg}(t_o)$; therefore it does not depend on the transaction intensity). For the less liquid stocks HON and NKE in Figures~\ref{fig:FigureVola3b}--\ref{fig:FigureVola4b} we chose $M=300$, $m = \lfloor 300\cdot \set{\text{\# of trades}} /T \rfloor$ (resulting in $m=53$ for HON; $m=55$ for  NKE) and $H=15$. The weighting functions $\fK (x)$, $K(x)$, and $k(x)$ applied here are Epanechnikov kernels $\frac{3}{4} \bra{ 1-x^2} I_{ \set{ \abs{x} \leq 1}}$ and the weighting function $g(x)$ in the pre-averaging steps is $g(x) = x(1-x)I_{\set{0 \leq x \leq 1}}$.

The bandwidth for the green estimator in Figures~\ref{fig:FigureVola1b}--6 is chosen adaptively as described in Dahlhaus and Neddermeyer (2013). This (recursive) estimator is equivalent to a kernel estimator with a one-sided kernel which results in an additional bias. This bias has been corrected in the above plots by a time-shift.

For the simulation in Figures~5--6 we have used additive noise with rounding, namely
$$
Y_{t_i} = \ \log\bra{ \lfloor 100\cdot(\exp(X_{t_i})+\varepsilon_{i} ) \rfloor /100 }
$$
with $\varepsilon_{i} \sim \cN(0,0.001^2)$. The tick-time volatility is chosen to be $\sigma^2(t) = \exp(-18)$ in Figure~5 and $\sigma^2(t) = \exp(-18 + \cos(10\pi t/T))$ in Figure~6 to mimic the shape of the empirical volatility per tick as in MSFT and NKE respectively. Furthermore, the transaction arrivals used here are taken from a real stock - CSCO on April 1, 2014, in order to mimic a real trading intensity.

\subsection*{A.2  Proof of the Volatility Decomposition}   
\small

\textbf{Proof of Proposition~\ref{prop_2_2}.} For $\delta >0$ and $t\in (0,T)$, we obtain with
$$t_{i_l} := \inf\{t_i : t_i > t\} \quad \text{and} \quad t_{i_u} := \sup\{t_i : t_i  \leq t+\delta\} $$
\begin{align*}
\bE &\brac{ \bra{X_{t+\delta} - X_t}^2 \con \cF_t} \ = \ \bE\brac{  \bra{\sum_{t<t_i \leq t+\delta} \sigma_{t_i} U_i}^2\con \cF_t} \\
&\quad  \stackrel {(*)} {=}  \ \bE\brac{  \bra{\sum_{t<t_i < t_{i_u}} \sigma_{t_i} U_i}^2 \con \cF_t }  + \bE\brac{  2\bra{\sum_{t<t_i < t_{i_u}} \sigma_{t_i} U_i}  \bE\brac{\sigma_{t_{i_u}} U_{i_u} \con \cF_{t_{i_u}-}} \con \cF_t}  \\
&\qquad + \ \bE\brac{ \  \bE\brac{ \sigma^2_{t_{i_u}} U^2_{i_u} \con \cF_{t_{i_u}-}} \con \cF_t} \\
&\quad = \ \bE\brac{  \bra{\sum_{t<t_i < t_{i_u}} \sigma_{t_i} U_i}^2 \con \cF_t }+ \bE\brac{\sigma^2_{t_{i_u}} \con \cF_t} \\
&\quad = \ \cdots = \ \bE\brac{ \sigma^2_{t_{i_u}} + \sigma^2_{t_{i_u -1}} + \ldots + \sigma^2_{t_{i_l}}  \con \cF_t}\\
&\quad = \ \bE\brac{ \int^{t+\delta}_t \sigma^2_s \ dN_s \con \cF_t} \ = \ \bE\brac{ \int^{t+\delta}_t \sigma^2_s \cdot \lambda_s \ ds \con \cF_t}
\end{align*}
which implies the assertion due to the dominated convergence theorem as the processes $\sigma_t$ and $\lambda_t$ are continuous over $[0,T]$. The last equality holds by the martingale property of $\int^t_0 \sigma_s dM_s$, since $M_t := N_t -\int^t_0 \lambda_s ds$ is a martingale and therefore the stochastic integral with $M_t$ as integrator is also a martingale. In $(*)$ we can take $\sigma_{t_{i_u}}$ out of the conditional expectation, since we have assumed that $\sigma_t$ is $\cF_t$-predictable, so $\sigma_{t_i}$ is $\cF_{t_i-}\,$-measurable. 
\qed
\normalsize

\subsection*{A.3  Proofs for the Rescaled Model}   

\small
In what follows, some steps of the proofs are related to martingale theory. In fact, the same results can be obtained without using it, however the calculations could be long and cumbersome. Another benefit of using this theory in our framework is that we can extend the proofs more easily to a more general case of stochastic intensity models where martingale dynamics are needed. We define $M_{t,T}$:=$N_{t,T}-\int^t_0 \lambda(s/T)ds$, which is clearly a martingale.  Then, it is possible to define a stochastic integral $\int^t_0 c_{s,T} \ dM_{s,T}$, where $c_{t,T}$ is a predictable process. We can show that
\begin{equation*} 
\bE \brac{ \int^t_0 c_{s,T} \ dM_{s,T} }^2 = \ \bE\brac{ \int^t_0 c_{s,T} \lambda(s/T) \ ds} \tag{A.1}
\end{equation*}
(e.g. Kuo 2006, chap.~6). Throughout the proofs, $C$ is used as a generic constant.\\

\noindent \textbf{Proof of Theorem~\ref{Thm:IntensityEstimator}.} 
By the isometry (A.1) we have 
\begin{align*}
\bE \brac{\widehat{\lambda}(u_o) - \frac{1}{\fb T} \int^T_0 \fK\bra{ \frac{t-u_oT }{\fb T}} \lambda\bra{ \frac{t}{T}} dt }^{2}  &= \frac{1}{\fb^2T^2}  \bE \brac{ \int^T_0 \fK\bra{ \frac{t-u_oT}{\fb T}} d M_{t,T} }^{2} \\
& = \frac{1}{\fb^2T^2} \int^T_0 \fK^{2}\bra{ \frac{t-u_oT}{\fb T}} \lambda\left( \frac{t}{T}\right)dt \ \to \ 0  
\end{align*}
for $\fb \leq u_o \leq 1-\fb$ with $0< \fb \leq1/2$. A usual bias calculation gives
\begin{align*}
&\frac{1}{\fb T} \int^T_0 \fK\bra{ \frac{t-u_oT }{\fb T}} \lambda\bra{ \frac{t}{T} }dt - \lambda(u_o)  = \int_{\bR} \fK(x) \set{\lambda(u_o+x\fb)- \lambda(u_o)} dx,
\end{align*}
which implies the consistency of the intensity estimate, since $\lambda(\cdot)$ is bounded continuous. To show the asymptotic normality we use the limit theorem for triangular arrays. First we divide $[0,T]$ into $M_T$ equidistant subintervals of a fixed length $\triangle$, i.e.~$M_T= \lfloor T/ \triangle\rfloor \to \infty$, as $T\to \infty$. For $j=0,1,...,M_T$ we define $\triangle_j:=j\triangle$ and rewrite 
\begin{align*}
\widehat{\lambda}(u_o) &= \sum^{M_T}_{i=1} \int^{\triangle_i}_{\triangle_{i-1} } \frac{1}{\fb T} \fK\bra{ \frac{t-u_oT}{\fb T}} dN_{t,T} +  \int^T_{\triangle_{M_T} }\frac{1}{\fb T} \fK\bra{\frac{t-u_oT}{\fb T}} dN_{t,T},
\end{align*}
which gives
\begin{align*} 
\sqrt{\fb T} \left( \widehat{\lambda}(u_o) - \frac{1}{\fb T} \int^{T}_{0} \fK\bra{ \frac{t-u_oT}{\fb T}} \lambda\bra{ \frac{t}{T} } dt\right) &= \sum^{M_T}_{i=1} \int^{\triangle_i}_{\triangle_{i-1} } \frac{1}{\sqrt{\fb T}} \fK\bra{ \frac{t-u_oT}{\fb T}} dM_{t,T}  +  R \\
 &=: \sum^{M_T}_{i=1}  Z_{i,T}  +  R. \tag{A.2} 
\end{align*}

\noindent Owing to independent increments of $N_{t,T}$, the random variables $Z_{i,T}$, $i=1,...,M_T,$ are independent with $\bE[Z_{i,T}]=0$. We can show that (i)
\begin{align*}
&V^2_{M_T} := \sum^{M_T}_{i=1} \bE[Z^2_{i,T}] \to \ \lambda(u_o) \int_{\bR} \fK^{2}(x) dx \ =: \ V^2;
\end{align*}
and (ii) (Lindeberg's condition) for all $\varepsilon >0$,
\begin{align*}
&\frac{1}{V^2_{M_T}} \sum^{M_T}_{i=1} \bE\left[ Z^2_{i,T} \ I_{\set{|Z_{i,T}| \geq \varepsilon\cdot V_{M_T}}}\right] =\frac{1}{V^2_{M_T}}  \sum^{M_T}_{i=1} \bE\left[ Z^2_{i,T} \ I_{ \set{ \abs{\int^{\triangle_i}_{\triangle_{i-1} }  \fK\bra{ \frac{t-u_oT}{\fb T}} dM_{t,T} }\geq  \sqrt{bT} \cdot\varepsilon\cdot V_{M_T}}}\right] \to \ 0
\end{align*}
by the dominated convergence theorem, since $\bE \brac{ Z^2_{i,T} } < \infty$. Thus $\sum^{M_T}_{i=1} Z_{i,T}  \xrightarrow{\cD} \cN(0,V^2)$ which implies that (A.2) $\stackrel{\cD}{\to} \cN(0,V^2)$, since the rest term $R$ is asymptotically negligible. Similarly, 
\begin{align*}
\bE \brac{\widehat{\lambda}(u_o)}  -\lambda(u_o)
=& \ \frac{1}{2} \lambda^{(2)}(u_o) \int_{\bR} \fb^2 x^2\fK (x)dx \cdot I_{\set{m' =2}} + O\bra{ \int_\bR \fK (x) |x\fb |^{m'+\gamma'} dx},
\end{align*}
since $\lambda(\cdot)$ lies in $\cC^{m',\gamma'}$. The assertion is then verified by the condition $\fb^{2(m'+\gamma')+1}T = o(1)$. 
\qed 


\noindent \textbf{Proof of Theorem~\ref{Thm:ClockTimeVolatilityEstimator}.} 
First we will focus on the asymptotic normality of the statistic
\begin{align*}
\widehat{\sigma}^2_{clock,pre}(u_o) :=& \frac{1}{b} \frac{1}{g_2}\sum^{M_T -1}_{i=0} K\bra{ \frac{t_{iH} - u_oT}{bT}} \bra{ \overline{\triangle Y}_{t_{iH},T} }^2  \\
&\quad -  \frac{1}{2bH} \frac{\sum^{H-1}_{l=1} h^2(l/H)}{ g_2} \sum^{N_T}_{i=1} K\bra{ \frac{t_i - u_oT}{bT}} \bra{Y_{t_i,T} - Y_{t_{i-1},T}}^2
\end{align*}
with $M_{T}:= \lfloor N_{T}/H \rfloor$ being the random number of blocks. The first term consists of non-overlapping blocks of data so that, conditionally on $N_{\cdot,T}$, a central limit theorem for independent triangular arrays can be applied. The second term of $\widehat{\sigma}^2_{clock,pre}(\cdot)$ remains the same as that of $\widehat{\sigma}^2_{clock,pavg}(\cdot)$ and plays no role in the limit distribution (it corrects the bias caused by the additive microstructure noise). In the end of the proof, we will show that the distinction between the two estimators $\widehat{\sigma}^2_{clock,pre}(\cdot)$ and $\widehat{\sigma}^2_{clock,pavg}(\cdot)$ is asymptotically negligible so that both have the same limit distribution.

We divide the proof into the following parts. Parts $(I)-(III)$ concern the limit distribution, where the rate of convergence will be balanced by the choice of the block size $H$. Parts $(IV)-(V)$ are related to biases; especially $(V)$ deals with the bias caused by the additional noise. 
\begin{align*}
&(I) \ \sqrt{\frac{bT}{H}} \Biggl( \frac{1}{b} \frac{1}{g_2} \sum^{M_T-1}_{i=0} K\bra{  \frac{t_{iH} - u_oT}{bT}} \Biggl( \bra{\ot{X}_{t_{iH},T}}^2 - \bE \brac{  \bra{\ot{X}_{t_{iH},T}}^2 \con N_{\cdot,T} }  \Biggr) \Biggr) \xrightarrow{\cD} \cN(0, \eta^2_A),\\
&(II) \ \sqrt{bH} \Biggl( \frac{1}{b} \frac{2}{g_2} \sum^{M_T-1}_{i=0} K\bra{  \frac{t_{iH} - u_oT}{bT}} \bra{\ot{X}_{t_{iH},T}} \bra{\ot{\ep}_{iH}} \Biggr)  \xrightarrow{\cD} \cN(0, \eta^2_B),  \quad \text{ and} \\
&(III) \ \sqrt{\frac{bH^3}{T}} \Biggl( \frac{1}{b} \frac{1}{g_2} \sum^{M_T-1}_{i=0} K\bra{  \frac{t_{iH} - u_oT}{bT}} \Biggl( \bra{\ot{\ep}_{iH}}^2 - \bE\brac{\bra{\ot{\ep}_{iH}}^2}   \Biggr) \Biggr) \xrightarrow{\cD} \cN(0, \eta^2_C).
\end{align*}
Furthermore, under the condition $H=\delta\cdot T^{1/2}$ for $\delta\in (0,\infty)$ and the condition $b^{2\alpha+1}T^{1/2} \to 0$ with $\alpha = \min\set{m+\gamma, m'+\gamma'}$, the biases are negligible in the limit, i.e. 
\begin{align*}
&(IV) \ \sqrt{bT^{1/2}} \Biggl( \frac{1}{b} \frac{1}{g_2} \sum^{M_T-1}_{i=0} K\bra{  \frac{t_{iH} - u_oT}{bT}} \bE \brac{  \bra{\ot{X}_{t_{iH},T}}^2 \con N_{\cdot,T}} -\sigma^2(u_o)\lambda(u_o) - BIAS \Biggr) = \ o_p(1), \\
&(V) \ \sqrt{bT^{1/2}} \Biggl(  \frac{1}{b} \frac{1}{g_2} \sum^{M_T-1}_{i=0} K\bra{  \frac{t_{iH} - u_oT}{bT}} \bE \brac{  \bra{\ot{\ep}_{iH}}^2 } \\
&\qquad \qquad - \frac{1}{2bH} \frac{\sum^{H-1}_{l=1} h^2(l/H)}{g_2} \sum^{N_T}_{i=1} K\bra{ \frac{t_i - u_oT}{bT}} \bra{Y_{t_i,T} - Y_{t_{i-1},T}}^2  \Biggr)  = \  o_p(1).
\end{align*}

We set the left-hand side of $(I)$ to be $\sum^{M_T-1}_{i=0} A_{i,T}$ where
\begin{align*}
&A_{i,T} = \sqrt{\frac{bT}{H}} \frac{1}{b} \frac{1}{g_2} K\bra{  \frac{t_{iH} - u_oT}{bT}} \Biggl( \bra{\ot{X}_{t_{iH},T}}^2 - \bE \brac{  \bra{\ot{X}_{t_{iH},T}}^2 \con N_{\cdot,T} } \Biggr). 
\end{align*}
We see that conditional on $N_{\cdot, T}$, $\set{A_{i,T}}_{i=0,...,M_T-1}$ is an independent sequence, so it is sufficient to show that (Hall and Heyde 1980, cor.~3.1)
\begin{align*}
(a) \ \sum^{M_T-1}_{i=0} \bE\brac{ A^2_{i,T} \con N_{\cdot,T}} \xrightarrow{\bP} \eta^2_{A} \ \text{ and }  \ (b) \ \sum^{M_T-1}_{i=0} \bE\brac{ A^4_{i,T} \con N_{\cdot,T}} \xrightarrow{\bP} 0.
\end{align*}
\noindent Corresponding to independent increments of $X_{\cdot,T}$, it gives
\begin{align*}
&\sum^{M_T-1}_{i=0} \bE\brac{ A^2_{i,T} \con N_{\cdot,T}} \\
&\quad = \sum^{M_T-1}_{i=0}  \frac{T}{bH}  \frac{1}{g^2_2} K^2 \bra{  \frac{t_{iH} - u_oT}{bT}} \set{ \bE \brac{\bra{\ot{X}_{t_{iH},T}}^4 \con N_{\cdot,T} }- \bra{\bE \brac{  \bra{\ot{X}_{t_{iH},T}}^2 \con N_{\cdot,T} }}^2 } \\
&\quad = \sum^{M_T-1}_{i=0}  \frac{2}{H}\frac{1}{bT} \frac{1}{g^2_2} K^2 \bra{  \frac{t_{iH} - u_oT}{bT}}  \set{ \sum^{H-1}_{l=1} g^4\bra{\frac{l}{H}} \sigma^4\bra{ \frac{t_{iH+l}}{T}} + \bra{ \sum^{H-1}_{l=1} g^2\bra{\frac{l}{H}} \sigma^2\bra{ \frac{t_{iH+l}}{T}}}^2 } \tag{see (A.3) and (A.4) below}\\
&\quad = \sum^{N_T}_{i=1}   \frac{2}{H^2}\frac{1}{bT} \frac{1}{g^2_2} K^2 \bra{  \frac{t_i - u_oT}{bT}} \sigma^4\bra{\frac{t_i}{T}} \set{ \sum^{H-1}_{l=1} g^4\bra{\frac{l}{H}}  +\bra{ \sum^{H-1}_{l=1}  g^2\bra{\frac{l}{H}} }^2} \ + \ o_p(1) \tag{see (A.5) below}\\
&\quad = \frac{2}{H^2} \frac{(\sum g^2(l/H))^2}{g^2_2} \sum^{N_T}_{i=1} \frac{1}{bT} K^2\bra{\frac{t_i -u_oT}{bT}} \sigma^4\bra{\frac{t_i}{T}} \ + \ o_p(1) \\ 
&\quad \xrightarrow{\bP} \ 2\sigma^4(u_o)\lambda(u_o) \int_\bR K^2(x)dx.
\end{align*}
\noindent By independent increments of $X_{\cdot,T}$,
\begin{align*}
\bE\brac{ (\ot{X}_{t_{iH},T} )^4 \con N_{\cdot,T}}  &= \ 2 \sum^{H-1}_{l =1} g^4\bra{\frac{l}{H}}  \sigma^4 \bra{\frac{t_{iH+l}}{T}} \frac{1}{T^2} + 3 \bra{ \sum^{H-1}_{l =1} g^2\bra{\frac{l}{H}}  \sigma^2 \bra{\frac{t_{iH+l}}{T}} \frac{1}{T} }^2;  \tag{A.3}
\end{align*}
\begin{align*}
&\bra{\bE\brac{ (\ot{X}_{t_{iH},T} )^2 \con N_{\cdot,T}} }^2 = \bra{ \sum^{H-1}_{l =1} g^2\bra{\frac{l}{H}}  \sigma^2 \bra{\frac{t_{iH+l}}{T}} \frac{1}{T} }^2. \tag{A.4}
\end{align*}
\noindent As $K$ has bounded first derivatives and $\sigma(\cdot) \in \cC^{m,\gamma}$,
\begin{align*}
 H&\cdot \sum^{M_T-1}_{i=0}  \frac{1}{H}\frac{1}{bT} \frac{1}{g^2_2} K^2 \bra{  \frac{t_{iH} - u_oT}{bT}} \sum^{H-1}_{l=1} g^4\bra{\frac{l}{H}} \sigma^4\bra{ \frac{t_{iH+l}}{T}} \\
&- \sum^{N_T}_{j=1}  \frac{1}{H}\frac{1}{bT} \frac{1}{g^2_2} K^2 \bra{  \frac{t_j - u_oT}{bT}} \sum^{H-1}_{l=1} g^4\bra{\frac{l}{H}} \sigma^4\bra{ \frac{t_{\ubra{j}+l}}{T}} \\
&+ \sum^{N_T}_{j=1}  \frac{1}{H}\frac{1}{bT} \frac{1}{g^2_2} K^2 \bra{  \frac{t_j - u_oT}{bT}} \sum^{H-1}_{l=1} g^4\bra{\frac{l}{H}} \sigma^4\bra{ \frac{t_{\ubra{j} +l}}{T}} \\
&- \sum^{N_T}_{j=1}  \frac{1}{H}\frac{1}{bT} \frac{1}{g^2_2} K^2 \bra{  \frac{t_j - u_oT}{bT}} \sum^{H-1}_{l=1} g^4\bra{\frac{l}{H}} \sigma^4\bra{ \frac{t_j}{T}}\\
=& \ O\bra{\frac{H}{bT}} + O_p\bra{ \abs{ \frac{H}{T}}^\gamma} \ = \ o_p(1) \tag{A.5},
\end{align*}
since $H/bT \to 0$. We now turn to condition $(b)$:
\begin{align*}
&\sum^{M_T-1}_{i=0} \bE\brac{ A^4_{i,T} \con N_{\cdot,T}} \\
&\quad \leq C\cdot \frac{T^2}{b^2H^2}  \sum^{M_T-1}_{i=0}  K^4 \bra{  \frac{t_{iH} - u_oT}{bT}}  \bE \brac{ \set{ \bra{\ot{X}_{t_{iH},T}}^2 - \bE \brac{  \bra{\ot{X}_{t_{iH},T}}^2 \con N_{\cdot,T} }}^4 \con N_{\cdot,T}} \\
&\quad \leq C\cdot \frac{T^2}{b^2H^2}  \sum^{M_T-1}_{i=0}  K^4 \bra{  \frac{t_{iH} - u_oT}{bT}} \bra{\sum^{H-1}_{l=1}g^2\bra{\frac{l}{H}}}^4 \frac{1}{T^4} \ = \ o_p(1),
\end{align*}
since $H/bT \to 0$ and $\bE \brac{  \bra{\ot{X}_{t_{iH},T}}^8 \con N_{\cdot,T} }  \leq C\cdot \set{\sum^{H-1}_{l=1} g^2\bra{\frac{l}{H}} \sigma^2\bra{\frac{t_{iH+l}}{T}} }^4\frac{1}{T^4}$ holds. In order to derive the limit distribution in $(II)$ we set
$$
\sqrt{bH} \Biggl( \frac{1}{b} \frac{2}{g_2} \sum^{M_T-1}_{i=0} K\bra{  \frac{t_{iH} - u_oT}{bT}} \bra{\ot{X}_{t_{iH},T}} \bra{\ot{\ep}_{iH}}  \Biggr)  =:  \sum^{M_T-1}_{i=0} B_{i,T}.
$$
Likewise we show that (conditional on $N_{\cdot,T}$)
$$
(c) \ \sum^{M_T-1}_{i=0} \bE\brac{ B^2_{i,T} \con N_{\cdot,T}} \xrightarrow{\bP} \eta^2_{B}  \ \text{ and } \
(d) \ \sum^{M_T-1}_{i=0} \bE\brac{ B^4_{i,T} \con N_{\cdot,T}} \xrightarrow{\bP} 0.
$$
Since 
$$
\bE \brac{ \bra{\ot{\ep}_{iH}}^2}  = \ \bE \brac{ \bra{ \sum^{H-1}_{l=1} h\bra{\frac{l}{H}}\ep_{iH+l} }^2}  = \ \omega^2 \sum^{H-1}_{l=1} h^2\bra{\frac{l}{H}},
$$
it implies that (see (A.4))
\begin{align*}
\bE\brac{ B^2_{i,T} \con N_{\cdot,T}}  &=  \frac{H}{b} \frac{4}{g^2_2}  K^2\bra{  \frac{t_{iH} - u_oT}{bT}}  \bE\brac{\bra{\ot{X}_{t_{iH},T}}^2 \con N_{\cdot,T}}  \cdot \bE\brac{ \bra{\ot{\ep}_{iH}}^2 \con N_{\cdot,T}}\\
& = \frac{H}{b} \frac{4}{g^2_2}  K^2\bra{  \frac{t_{iH} - u_oT}{bT}}  \bra{ \sum^{H-1}_{l=1}g^2\bra{\frac{l}{H}} \sigma^2\bra{\frac{t_{iH+l}}{T}} \frac{1}{T}} \bra{ \omega^2 \sum^{H-1}_{l=1} h^2\bra{\frac{l}{H}} } .
\end{align*}
Hence (similar to (A.5))
\begin{align*}
\sum^{M_T-1}_{i=0} &\bE\brac{ B^2_{i,T} \con N_{\cdot,T}}  \ \xrightarrow{\bP} \ 4\omega^2 \sigma^2(u_o)\lambda(u_o) \frac{g'_2}{g_2}\int_\bR K^2(x)dx, 
\end{align*}
since $ \sum^{H-1}_{l=1} h^2\bra{\frac{l}{H}}  = \sum^{H-1}_{l=1} \set{ g^{(1)}\bra{\frac{l}{H}}\frac{1}{H} + o\bra{\frac{1}{H}} }^2.$ Similarly,
\begin{align*}
&\sum^{M_T-1}_{i=0} \bE\brac{ B^4_{i,T} \con N_{\cdot,T}} = \sum^{M_T-1}_{i=0} \frac{H^2}{b^2} \frac{16}{g^4_2}  K^4\bra{  \frac{t_{iH} - u_oT}{bT}} \bE\brac{\bra{\ot{X}_{t_{iH},T}}^4 \con N_{\cdot,T}}  \cdot \bE\brac{ \bra{\ot{\ep}_{iH}}^4 \con N_{\cdot,T}}\\
&\quad \leq \ C \cdot \sum^{M_T-1}_{i=0} \frac{H^2}{b^2}  K^4\bra{  \frac{t_{iH} - u_oT}{bT}} \bra{ \sum^{H-1}_{l=1}g^2\bra{\frac{l}{H}}\frac{1}{T}}^2 \bra{ \sum^{H-1}_{l=1} h^2\bra{\frac{l}{H}} }^2 \ \xrightarrow{\bP} \ 0. 
\end{align*}
We proceed analogously to show $(III)$. We denote its left-hand side by $\sum^{M_T-1}_{i=0} C_{i,T}$ and show that 
$$
(e) \sum^{M_T-1}_{i=0} \bE\brac{ C^2_{i,T} \con N_{\cdot,T}} \xrightarrow{\bP} \eta^2_{C} \ \text{ and } \ (f) \sum^{M_T-1}_{i=0} \bE\brac{ C^4_{i,T} \con N_{\cdot,T}} \xrightarrow{\bP} 0.
$$

\noindent On account of the assumption $\bE\brac{\ep^4_i} = \theta\omega^4$, $\theta \in \bR^+$, we have
\begin{align*}
\bE\brac{(\ot{\ep}_{iH})^4} 
=& \ (\theta -1) \sum^{H-1}_{l=1} h^4\bra{\frac{l}{H}} \omega^4  + 3 \bra{\sum^{H-1}_{l=1} h^2\bra{\frac{l}{H}} }^2 \omega^4,
\end{align*}
thus
\begin{align*}
&\sum^{M_T-1}_{i=0} \bE\brac{ C^2_{i,T} \con N_{\cdot,T}}  
= \sum^{M_T-1}_{i=0} \frac{H^3}{bT} \frac{\omega^4}{g^2_2} K^2\bra{  \frac{t_{iH} - u_oT}{bT}}  \set{ (\theta -1) \sum^{H-1}_{l=1} h^4\bra{\frac{l}{H}}   + 2 \bra{\sum^{H-1}_{l=1} h^2\bra{\frac{l}{H}} }^2  } \\
&\quad = \ \frac{1}{H} \cdot H\sum^{M_T-1}_{i=0} \frac{1}{bT}K^2\bra{  \frac{t_{iH} - u_oT}{bT}} \frac{(\theta-1)\omega^4}{g^2_2} \set{ \sum^{H-1}_{l=1} \bra{g'\bra{\frac{l}{H}}}^4\frac{1}{H} + o\bra{\frac{1}{H}} }\\
&\qquad + \frac{2\omega^4}{g^2_2} \cdot H\sum^{M_T-1}_{i=0} \frac{1}{bT}K^2\bra{  \frac{t_{iH} - u_oT}{bT}}  \set{ \sum^{H-1}_{l=1} \bra{g'\bra{\frac{l}{H}}}^2\frac{1}{H} + o\bra{\frac{1}{H}} }^2\\
&\quad \xrightarrow{\bP} \ 2\omega^4 \lambda(u_o) (g'_2/g_2)^2\int_\bR K^2(x)dx,
\end{align*}
which leads to $(e)$. In fact, $(f)$ can be established in the same manner as before, therefore omitted. To build the joint distribution of $(I), \ (II)$ and $(III)$ it suffices to show that 
$$
\sum^{M_T-1}_{i=0} a\sqrt{\delta}A_{i,T}+ b\frac{1}{\sqrt{\delta}} B_{i,T} + c \frac{1}{\sqrt{\delta^3}} C_{i,T} \ \xrightarrow{\cD} \ a\sqrt{\delta}A+ b\frac{1}{\sqrt{\delta}} B + c \frac{1}{\sqrt{\delta^3}} C
$$
by Cramer-Wold's theorem, for all $ a,b,c \in \bR$, where $A, B$ and $C$ are the limits of $(I), \ (II)$ and $(III)$ respectively. Indeed, this joint limit is a direct consequence of $a) - f)$ where now the pre-averaging block size $H$ is chosen to equal $\delta \cdot T^{1/2}$. To sum up, we have shown the first main part of the limit distribution, i.e. 
\begin{align*}
&\sqrt{bT^{1/2}} \Biggl( \sum^{M_T-1}_{i=0} \frac{1}{b} \frac{1}{g_2} K\bra{  \frac{t_{iH} - u_oT}{bT}}  \set{ \bra{\ot{Y}_{t_{iH},T}}^2 - \bE \brac{  \bra{\ot{X}_{t_{iH},T}}^2 \con N_{\cdot,T} } -\bE\brac{\bra{\ot{\ep}_{iH}}^2} } \Biggr) \\
&\qquad \xrightarrow{\cD} \ \cN\bra{ 0, \ \delta\eta^2_A +\frac{1}{\delta}\eta^2_B + \frac{1}{\delta^3}\eta^2_C },
\end{align*}
since $\bra{\ot{Y}_{t_{iH},T}} = \bra{\ot{X}_{t_{iH},T}} + \bra{\ot{\ep}_{t_{iH},T}}$. 

We carry out the proof by showing the asymptotic biases in $(IV)$ and $(V)$. In order to derive $(IV)$ we see that
\begin{align*}  
&\sqrt{bT^{1/2}} \Biggl( \frac{1}{b} \frac{1}{g_2} \sum^{M_T-1}_{i=0} K\bra{  \frac{t_{iH} - u_oT}{bT}} \bE \brac{  \bra{\ot{X}_{t_{iH},T}}^2 \con N_{\cdot,T}}  -\sigma^2(u_o)\lambda(u_o) - BIAS \Biggr) \\
&= \sqrt{bT^{1/2}}  \frac{1}{H} \frac{1}{g_2} \cdot\Biggl( H \sum^{M_T-1}_{i=0} \frac{1}{bT} K\bra{  \frac{t_{iH} - u_oT}{bT}} \sum^{H-1}_{l=1} g^2\bra{\frac{l}{H}} \sigma^2\bra{\frac{t_{iH+l}}{T}} \\
&\qquad - \sum^{N_T}_{j=1} \frac{1}{bT} K\bra{  \frac{t_j - u_oT}{bT}} \sum^{H-1}_{l=1} g^2\bra{\frac{l}{H}} \sigma^2\bra{\frac{t_{\ubra{j}+l}}{T}} \Biggr) \notag\\
& \quad + \sqrt{bT^{1/2}}  \frac{1}{H} \frac{1}{g_2} \sum^{N_T}_{j=1} \frac{1}{bT} K\bra{  \frac{t_j - u_oT}{bT}}  \Biggl(  \sum^{H-1}_{l=1} g^2\bra{\frac{l}{H}} \sigma^2\bra{\frac{t_{\ubra{j}+l}}{T}}  -   \sigma^2\bra{\frac{t_j}{T}}\sum^{H-1}_{l=1} g^2\bra{\frac{l}{H}}  \Biggr) \\
& \quad + \sqrt{bT^{1/2}}  \Biggl( \frac{1}{H} \frac{1}{g_2} \sum^{N_T}_{j=1} \frac{1}{bT} K\bra{  \frac{t_j - u_oT}{bT}}  \sigma^2\bra{\frac{t_j}{T}}\sum^{H-1}_{l=1} g^2\bra{\frac{l}{H}} \\
&\qquad - \sigma^2(u_o)\lambda(u_o) - \frac{1}{2}\bra{ \sigma^2(u_o)\lambda(u_o)}^{(2)} b^2 \int_\bR x^2K(x)dx I_{\set{m=m'=2}} \Biggr) \\
& =  O\bra{ \sqrt{bT^{1/2}}\cdot \frac{H}{bT}} + O_p\bra{\sqrt{bT^{1/2}}\abs{\frac{H}{T}}^\gamma}  + (iv) \ = \ o_p(1), \tag{A.6}
\end{align*} 
since $K$ has first bounded derivatives and $\sigma(\cdot) \in \cC^{m,\gamma}$ (see also (A.5)) and
\begin{align*}
&(iv)  = \sqrt{bT^{1/2}} \Biggl( \frac{1}{H} \frac{1}{g_2} \sum^{N_T}_{j=1} \frac{1}{bT} K\bra{  \frac{t_j - u_oT}{bT}}  \sigma^2\bra{\frac{t_j}{T}}\sum^{H-1}_{l=1} g^2\bra{\frac{l}{H}} \\
&\qquad \qquad- \bE\brac{ \frac{1 }{H} \frac{1}{g_2} \sum^{N_T}_{j=1} \frac{1}{bT} K\bra{  \frac{t_j - u_oT}{bT}}  \sigma^2\bra{\frac{t_j}{T}}\sum^{H-1}_{l=1} g^2\bra{\frac{l}{H}} } \Biggr)\\
&\quad \quad + \sqrt{bT^{1/2}} \Biggl( \bE\brac{ \frac{1}{H} \frac{1}{g_2} \sum^{N_T}_{j=1} \frac{1}{bT} K\bra{  \frac{t_j - u_oT}{bT}}  \sigma^2\bra{\frac{t_j}{T}}\sum^{H-1}_{l=1} g^2\bra{\frac{l}{H}}} \\
& \qquad \qquad - \set{ \sigma^2(u_o)\lambda(u_o) + \frac{1}{2}\bra{ \sigma^2(u_o)\lambda(u_o)}^{(2)}b^2 \int_\bR x^2K(x)dx\cdot I_{\set{m=m'=2}}} \Biggr)\\
&= \sqrt{bT^{1/2}}\bra{ \int^T_0 \frac{1}{bT} K\bra{\frac{t-u_oT}{bT}} \sigma^2\bra{\frac{t}{T}} dM_{t,T} }\frac{\frac{1}{H} \sum g^2(l/H)}{g_2} +  O\bra{\sqrt{bT^{1/2}} b^{\min(\gamma, \gamma')}}  = \ o_p(1).
\end{align*}
The last equality in (A.6) is satisfied by the bandwidth condition $b^{2\alpha+1} T^{1/2} \to 0$. Finally, we separate $(V)$ into two summands
\begin{align*}
&\sqrt{bT^{1/2}} \Biggl(  \frac{1}{b} \frac{1}{g_2} \sum^{M_T-1}_{i=0} K\bra{  \frac{t_{iH} - u_oT}{bT}} \bE \brac{  \bra{\ot{\ep}_{iH}}^2 }  - \frac{T}{H} \frac{ \lambda(u_0)}{g_2} \sum^{H-1}_{l=1}h^2\bra{\frac{l}{H}}\omega^2 \Biggr) \\
& + \sqrt{bT^{1/2}} \Biggl(  \frac{T}{H} \frac{ \lambda(u_0)}{g_2} \sum^{H-1}_{l=1}h^2\bra{\frac{l}{H}}\omega^2 - \frac{T}{2H} \frac{\sum^{H-1}_{l=1} h^2(l/H)}{g_2} \cdot \frac{1}{bT} \sum^{N_T}_{l=1} K\bra{ \frac{t_i - u_oT}{bT}} \bra{Y_{t_i,T} - Y_{t_{i-1},T}}^2 \Biggr) \\
&=: (v_1) + (v_2).
\end{align*}
Since $\bE \brac{  \bra{\ot{\ep}_{iH}}^2 } = \sum^{H-1}_{l=1}h^2\bra{\frac{l}{H}}\omega^2$, direct calculations yield
\begin{align*}
&\bE \brac{ (v_1)^2} \leq  C\cdot bT^{1/2} \bra{\sum^{H-1}_{l=1}h^2\bra{\frac{l}{H}}}^2  \bE\brac{ \frac{1}{b} \sum^{M_T-1}_{i=0} K\bra{  \frac{t_{iH} - u_oT}{bT}}  - \frac{T}{H} \lambda(u_o) }^2 \\
&\leq  C\cdot  bT^{1/2}\frac{T^2}{H^2} \bra{\sum^{H-1}_{l=1}h^2\bra{\frac{l}{H}}}^2\Biggl( \bE\brac{\frac{H}{bT} \sum^{M_T-1}_{i=0} K\bra{ \frac{t_{iH} - u_oT}{bT}}  - \bE\brac{\frac{H}{bT} \sum^{M_T-1}_{i=0} K\bra{  \frac{t_{iH} - u_oT}{bT}} } \  }^2 \\
& \qquad + \bra{ \bE\brac{\frac{H}{bT} \sum^{M_T-1}_{i=0} K\bra{  \frac{t_{iH} - u_oT}{bT}}}  - \lambda(u_o)}^2 \Biggr)\\
&=  O\bra{bT^{1/2} \frac{T^2}{H^2} \frac{1}{H^2} \frac{1}{bT}} + O\bra{ bT^{1/2} \frac{T^2}{H^2} \frac{1}{H^2}  b^{2\gamma'}}  = \ o(1).
\end{align*}
The last term $(v_2)$ is divided again into the sum of $(v_{2,1})$ and $(v_{2,2})$, where
\begin{align*}
&(v_{2,1}) := \ \sqrt{bT^{1/2}} \Biggl[  \frac{T}{H} \frac{ \lambda(u_0)}{g_2} \sum^{H-1}_{l=1}h^2\bra{\frac{l}{H}}\omega^2 - \frac{T}{2H} \frac{\sum^{H-1}_{l=1} h^2(l/H)}{g_2} \frac{1}{bT} \sum^{N_T}_{l=1} K\bra{ \frac{t_i - u_oT}{bT}} \bra{\ep_i - \ep_{i-1}}^2 \Biggr] 
\end{align*}
and
\begin{align*}
&(v_{2,2}):= -\sqrt{bT^{1/2}} \Biggl[ \frac{T}{2H} \frac{\sum^{H-1}_{l=1} h^2(l/H)}{g_2} \frac{1}{bT} \sum^{N_T}_{l=1} K\bra{ \frac{t_i - u_oT}{bT}} \times \\
&\qquad \qquad \times \bra{ \bra{X_{t_i,T}-X_{t_{i-1},T} }^2 + 2\bra{X_{t_i,T}-X_{t_{i-1},T} }\bra{\ep_i - \ep_{i-1}}}  \Biggr]. 
\end{align*}
It is easy to show that $\bE\brac{(v_{2,1})^2} = o(1)$ and $\bE\brac{ \ |(v_{2,2})| \ } = o(1)$ hold, therefore $(V)=o_p(1)$. In summary, we have shown the asymptotic normality for $\widehat{\sigma}^2_{clock,pre}(u_o)$.\\

To complete the proof, the difference between $\widehat{\sigma}^2_{clock,pavg}(u_o)$ and $\widehat{\sigma}^2_{clock,pre}(u_o)$ needs to be determined; in particular we show that 
\begin{equation*} 
\sqrt{bT^{1/2}} \set{ \widehat{\sigma}^2_{clock,pre}(u_o) - \widehat{\sigma}^2_{clock,pavg}(u_o)} \ = \ o_p(1).
\end{equation*}
Since the derivative of $K$ is bounded, it is enough to verify that (see also (A.5))
\begin{align*} 
&\sqrt{bT^{1/2}} \frac{1}{bH} \frac{1}{g_2}\sum^{N_T}_{i=1} K\bra{ \frac{t_i - u_oT}{bT}}  \set{ \bra{ \overline{\triangle Y}_{t_{\ubra{i}},T} }^2 - \bra{ \overline{\triangle Y}_{t_i,T} }^2 }  \\
&= \sqrt{bT^{1/2}} \frac{1}{bH} \frac{1}{g_2}\sum^{N_T}_{i=1} K\bra{ \frac{t_i - u_oT}{bT}} \Biggl( \ \set{ \bra{ \overline{\triangle X}_{t_{\ubra{i}},T} }^2 - \bra{ \overline{\triangle X}_{t_i,T} }^2 } +\set{ \bra{ \overline{\triangle \varepsilon}_{\ubra{i}} }^2 - \bra{ \overline{\triangle \varepsilon}_{i} }^2 } \\
&\quad + 2\set{ \bra{ \overline{\triangle X}_{t_{\ubra{i}},T} }\bra{ \overline{\triangle \varepsilon}_{\ubra{i}} } - \bra{ \overline{\triangle X}_{t_i,T} }\bra{ \overline{\triangle \varepsilon}_{i} } } \ \Biggr) \\
&=: \ (T_1) + (T_2) + (T_3) = \ o_p(1). \tag{A.7}
\end{align*}
We perform only the proof of the first term $(T_1)$. We will show below that $\bE\brac{(T_1)^2 \con N_{\cdot,T}} = o_p(1)$. 
\begin{align*}
&\bE \brac{(T_1)^2 \con N_{\cdot,T} }  = \ \frac{T^{1/2}}{bH^2} \frac{1}{g^2_2}\sum^{N_T}_{i,j=1} K\bra{ \frac{t_i - u_oT}{bT}} K\bra{ \frac{t_j - u_oT}{bT}} \times \\
&\quad  \times \Biggl( \bE \brac{ \bra{ \overline{\triangle X}_{t_{\ubra{i}},T} }^2 \bra{ \overline{\triangle X}_{t_{\ubra{j}},T} }^2 \con N_{\cdot,T}}  -  \bE \brac{ \bra{ \overline{\triangle X}_{t_{\ubra{i}},T} }^2 \bra{ \overline{\triangle X}_{t_j,T} }^2 \con N_{\cdot,T} } \\
&\qquad  \quad -  \bE \brac{ \bra{ \overline{\triangle X}_{t_i,T} }^2 \bra{ \overline{\triangle X}_{t_{\ubra{j}},T} }^2 \con N_{\cdot,T} } + \bE \brac{ \bra{ \overline{\triangle X}_{t_i,T} }^2 \bra{ \overline{\triangle X}_{t_j,T} }^2 \con N_{\cdot,T} } \Biggr) \\
&=:  \ (\triangle_1) - (\triangle_2) - (\triangle_3) +(\triangle_4) .
\end{align*}
The first term is split up into three small terms: 
\begin{align*}
&(\triangle_1) = \ \frac{T^{1/2}}{bH^2} \frac{1}{g^2_2}\sum_{i=j} K^2\bra{ \frac{t_i - u_oT}{bT}} \bE \brac{ \bra{ \overline{\triangle X}_{t_{\ubra{i}},T} }^2 \con N_{\cdot,T} } \\
& \ +  \frac{T^{1/2}}{bH^2} \frac{1}{g^2_2} \sum_{|i-j| \geq H} K\bra{ \frac{t_i - u_oT}{bT}} K\bra{ \frac{t_j - u_oT}{bT}}  \bE \brac{ \bra{ \overline{\triangle X}_{t_{\ubra{i}},T} }^2 \con N_{\cdot,T} } \cdot \bE \brac{ \bra{ \overline{\triangle X}_{t_{\ubra{j}},T} }^2 \con N_{\cdot,T} } \tag{since both terms are independent.} \\
& \ +  \frac{T^{1/2}}{bH^2} \frac{1}{g^2_2} \sum_{0<|i-j| < H} K\bra{ \frac{t_i - u_oT}{bT}} K\bra{ \frac{t_j - u_oT}{bT}}  \bE \brac{ \bra{ \overline{\triangle X}_{t_{\ubra{i}},T} }^2   \bra{ \overline{\triangle X}_{t_{\ubra{j}},T} }^2 \con N_{\cdot,T} } \\
&=: \ (\triangle_{1,1}) + (\triangle_{1,2}) + (\triangle_{1,3}).  
\end{align*}
Likewise, we expand $(\triangle_2) = (\triangle_{2,1}) + (\triangle_{2,2}) + (\triangle_{2,3})$. It is clear that $(\triangle_{1,1})$ and $(\triangle_{2,1})$ are of smaller order, therefore neglected. Since $(\triangle_{1,2})$ is equal to $(\triangle_{2,2})$, they cancel out. Lastly, 
\begingroup
\allowdisplaybreaks
\begin{align*}
&(\triangle_{1,3}) \ - \  (\triangle_{2,3})  \\ 
&\quad= \frac{T^{1/2}}{bH^2} \frac{1}{g^2_2} \sum_{0<|i-j| < H} K\bra{ \frac{t_i - u_oT}{bT}} K\bra{ \frac{t_j - u_oT}{bT}} \times \\
&\quad \qquad \times \bE \brac{ \bra{ \overline{\triangle X}_{t_{\ubra{i}},T} }^2 \set{\bra{ \overline{\triangle X}_{t_{\ubra{j}},T} }^2 - \bra{ \overline{\triangle X}_{t_j,T} }^2 }  \con N_{\cdot,T} } \\
&\quad= \frac{2T^{1/2}}{bH^2} \frac{1}{g^2_2} \sum^{N_T}_{i=1} \sum^{H-1}_{\alpha = 1} K\bra{ \frac{t_i - u_oT}{bT}} K\bra{ \frac{t_{i+\alpha} - u_oT}{bT}} \bE \Biggl[ \bra{ \overline{\triangle X}_{t_{\ubra{i}},T} }^2 \times\\
&\quad \qquad \times  \set{\bra{ \overline{\triangle X}_{t_{\ubra{(i+\alpha)}},T} }^2 - \bra{ \overline{\triangle X}_{t_{i+\alpha},T} }^2 }  \con N_{\cdot,T} \Biggr] \\
&\quad= \frac{2T^{1/2}}{bH^2} \frac{1}{g^2_2} \sum^{N_T}_{i=1}  K^2\bra{ \frac{t_i - u_oT}{bT}} \bE \Biggl[ \bra{ \overline{\triangle X}_{t_{\ubra{i}},T} }^2 \times\\
&\quad\qquad \times \underbrace{ \sum^{H-1}_{\alpha = 1} \set{\bra{ \overline{\triangle X}_{t_{\ubra{(i+\alpha)}},T} }^2 - \bra{ \overline{\triangle X}_{t_{i+\alpha},T} }^2 } }_{(\clubsuit)}  \con N_{\cdot,T} \Biggr]  + \ o(1)\\
&\quad= \frac{2T^{1/2}}{bH^2} \frac{1}{g^2_2} \sum^{N_T}_{i=1} K^2\bra{ \frac{t_i - u_oT}{bT}} o\bra{\frac{H^3}{T^2}} + o(1) \ = \ o_p(1), 
\end{align*}
\endgroup
since $g$ is differentiable, $g^{(1)}$ is Lipschitz's continuous, and 
\begin{align*}
(\clubsuit) &= \sum^{H-1}_{\alpha = 1}  \bra{ \overline{\triangle X}_{t_{\ubra{(i+\alpha)}},T}  - \overline{\triangle X}_{t_{i+\alpha},T} } \bra{ \overline{\triangle X}_{t_{\ubra{(i+\alpha)}},T}  + \overline{\triangle X}_{t_{i+\alpha},T} } \\
&= \sum^{H-1}_{\alpha = 1}  \sum^{H - (j \mod H)}_{l,l' = 1} \set{h\bra{ \frac{l+(j \mod H)}{H}} - h\bra{\frac{l}{H}}} \set{h\bra{ \frac{l'+(j \mod H)}{H}} - h\bra{\frac{l'}{H}}}  \times \\
&\qquad \times \bra{X_{t_{\alpha+l,T}} -X_{t_{\alpha+l-1,T}} } \bra{X_{t_{\alpha+l',T}} -X_{t_{\alpha+l'-1,T}} }.
\end{align*}
Thus, we get $(\triangle_1) - (\triangle_2) = o_p(1)$. Analogously,~$(\triangle_4) - (\triangle_3) = o_p(1)$ and hence $(T_1) = o_p(1)$. The rest terms $(T_2)$ and $(T_3)$ can be done in the same manner by employing $\overline{\triangle \varepsilon}_i = -\sum^{H-1}_{l=1} h(l/H) \varepsilon_{i+l}$ (also compare with $(II)$ and $(III)$). Therefore we conclude that (A.7) = $o_p(1)$.
\qed 

\medskip


For the bias derivation of the tick-time volatility estimate $\widehat{\sigma}^2_{pavg}(\cdot)$ we will need the following results - Lemma~1 and Corollary~1. In fact, these results are investigated under a general setting for point processes allowing for stochastic intensity. More precisely, given a filtered probability space $\bra{\Omega, \cF, \bra{\cF_{t,T}}_{t \in [0,T]}, \bP}$ a point process $N_{t,T}$ has an $\cF_{t,T}$-intensity $\lambda(t/T)$ if the conditions in Definition D7 Br\'emaud~(1981) holds. For a stopping time $\tau$ we define $\cF_{\tau,T}$ as consisting of sets $A\in \cF$ for which $A \cap \set{\tau \leq t} \in \cF_{t,T}$.

\noindent \textbf{Lemma 1} Suppose the intensity process $\lambda(u)$ is bounded continuous and bounded away from zero uniformly in $u\in[0,1]$, with probability 1. Then for $j \geq 0 $ and $0<l \leq 4$, it implies that 
\begin{align*}
\text{(i)}& \ \bE \brac{ \int^{t_{i+j}}_{t_i}\lambda(s/T) ds \con \cF_{t_i,T} } = j, \qquad \text{(ii)} \ \bE \brac{\bra{\int^{t_{i+j}}_{t_i}\lambda(s/T) ds}^2 \con \cF_{t_i,T}} = j^2+j \ \text{ and } \\
\text{(iii)}& \ \bE\brac{ (t_{i+j} - t_i)^l \con \cF_{t_i,T}} = O(j^l).
\end{align*}
\noindent \textbf{Proof.} Since $M_{t,T} = N_{t,T}-\int^t_0 \lambda(s/T)ds$ is a martingale and the arrival times $t_i$ are stopping times, we get
$$
j \ = \ \bE\brac{ N_{t_{i+j},T} - N_{t_i,T} \con \cF_{t_i,T}}  = \ \bE\brac{  \int^{t_{i+j}}_{t_i}\lambda(s/T) ds \con \cF_{t_i,T}}
$$
by the optional sampling theorem. Moreover, it is clear that $\tilde{M}_{t,T} := M^2_{t,T} - \int^t_0 \lambda(s/T) ds$ is another martingale, thus
\begin{align*}
0 =& \ \bE\brac{ \tilde{M}_{t_{i+j},T} - \tilde{M}_{t_i,T} \con \cF_{t_i,T}}  = \ \bE\brac{  \bra{ j - \int^{t_{i+j}}_{t_i}\lambda(l/T) dl}^2 \con \cF_{t_i,T}} -j, 
\end{align*}
i.e.~$\bE\brac{  \bra{\int^{t_{i+j}}_{t_i}\lambda(s/T) ds}^2 \con \cF_{t_i,T}} = j^2 + j$. Therefore we have shown (i) and (ii). Burkholder-Davis-Gundy's inequality (Jacod and Protter~2012, p.~39) yields
\begin{align*}
\bE\brac{  \abs{M_{t_{i+j},T} - M_{t_i,T}}^4 \con \cF_{t_i,T}} \leq&  \ C\cdot \bE\brac{  \bra{\int^{t_{i+j}}_{t_i}\lambda(s/T) ds}^2 \con \cF_{t_i,T}} = \ O(j^2).
\end{align*}
By applying H\"older's inequality, $\bE \brac{  \abs{M_{t_{i+j},T} - M_{t_i,T}}^3 \con \cF_{t_i,T}} = O(j^{3/2})$. Hence
\begin{align*}
\bE \brac{ (t_{i+j} - t_i)^4 \con \cF_{t_i,T}}  &\leq \ C\cdot \bE\brac{  \bra{\int^{t_{i+j}}_{t_i}\lambda(s/T) ds}^4 \con \cF_{t_i,T}}  \\
&= \ C\cdot \bE\brac{  \bra{M_{t_i,T} - M_{t_{i+j},T} \ + \ j}^4 \con \cF_{t_i,T}} = \ O(j^4).
\end{align*} 
The rest can be then easily justified by employing H\"older's inequality.
\qed \\


\noindent \textbf{Corollary 1} \ Let paths of $\lambda(\cdot) $ lie in $ \cC^{m',\gamma'}[0,1]$ and satisfy Assumption~\ref{assu_3_1}~iii). For $j \in \bN$, 
\begin{align*}
\text{(i)}& \ \bE\brac{ t_{i+j} - t_i \con \cF_{t_i,T}} = \frac{j}{\lambda(t_i/T)} \ - \ \frac{1}{2} \frac{ \lambda^{(1)}(t_i/T)}{\lambda^2(t_i/T)} \frac{j^2}{T}I_{\set{m'\neq 0}}  + O\bra{ \frac{j^{1+m'+\gamma'}}{T^{m'+\gamma'}}I_{\set{m'\neq 2}} } + O\bra{ \frac{j^3}{T^2}I_{\set{m'=2}} },\\
\text{(ii)}& \ \bE\brac{ (t_{i+j} - t_i)^2 \con \cF_{t_i,T}} = \frac{j^2 + j}{\lambda^2(t_i/T)}  + \ O\bra{ \frac{j^{2+\gamma'}}{T^{\gamma'}}I_{\set{m'=0}} + \frac{j^3}{T}I_{\set{m' \neq 0}} }.
\end{align*}

\noindent \textbf{Proof.}
By Lemma 1 we get
$$
\bE \brac{\lambda(t_i/T) (t_{i+j} - t_i) \con \cF_{t_i,T}} 
= \ \underbrace{\bE \brac{ \int^{t_{i+j}}_{t_i} \set{\lambda(t_i/T)- \lambda(l/T) }  dl \con \cF_{t_i,T}}}_{=:(\Lambda)} \ + \ j.
$$
Therefore, the first statement (i) is verified by considering the following cases: (for $\omega \in \Omega$)

\medskip

\noindent (a) for $\lambda(\cdot)(\omega) \in \cC^{0,\gamma'}$,
\begin{align*}
(\Lambda) \leq& \ C \cdot \brac{ \int^{t_{i+j}}_{t_i} \bra{ \frac{l-t_i}{T} }^{\gamma'}  dl \con \cF_{t_i,T}}  \ \leq \ C \cdot \frac{j^{1+\gamma'}}{T^{\gamma'} }; 
\end{align*}
(b) for $\lambda(\cdot)(\omega) \in \cC^{1,\gamma'}$,
\begin{align*}
(\Lambda) &= \ -\lambda^{(1)} \bra{ \frac{t_i}{T}} \bE\brac{ \int^{t_{i+j}}_{t_i} \bra{ \frac{l-t_i}{T} } dl \con \cF_{t_i,T}}   + O\bra{ \bE\brac{ \int^{t_{i+j}}_{t_i} \bra{ \frac{l-t_i}{T} }^{1+\gamma'} dl \con \cF_{t_i,T}}  }\\
&= -\frac{1}{2} \lambda^{(1)}\bra{ \frac{t_i}{T}} \bE\brac{ \frac{(t_{i+j} -t_i)^2}{ T} \con \cF_{t_i,T}} + O\bra{ \bE\brac{ \frac{(t_{i+j} -t_i)^{2+\gamma'}}{ T^{1+\gamma'}} \con \cF_{t_i,T}} } \\
&\stackrel{(*)}{=} \ -\frac{1}{2} \frac{ \lambda^{(1)}(t_i/T)}{\lambda^2(t_i/T)} \frac{j^2}{T}+ O\bra{ \frac{j^{2+\gamma'}}{T^{1+\gamma'}}};
\end{align*}
(c) for $\lambda(\cdot)(\omega) \in \cC^{2,\gamma'}$,
\begin{align*}
(\Lambda) &= \ -\lambda^{(1)} \bra{ \frac{t_i}{T}} \bE\brac{ \int^{t_{i+j}}_{t_i} \bra{ \frac{l-t_i}{T} } dl \con \cF_{t_i,T}} + O\bra{ \bE\brac{ \int^{t_{i+j}}_{t_i} \bra{ \frac{l-t_i}{T} }^2 dl \con \cF_{t_i,T}}  }\\
&= \ -\frac{1}{2} \lambda^{(1)}\bra{ \frac{t_i}{T}} \bE\brac{ \frac{(t_{i+j} -t_i)^2}{ T} \con \cF_{t_i,T}} + O\bra{ \bE\brac{ \frac{(t_{i+j} -t_i)^3}{ T^2} \con \cF_{t_i,T}} } \\
&\stackrel{(*)}{=} \ -\frac{1}{2} \frac{ \lambda^{(1)}(t_i/T)}{\lambda^2(t_i/T)} \frac{j^2}{T} + O\bra{ \frac{j^3}{T^2}}.
\end{align*}
To show $(*)$ in (b) and (c), we have to apply the result (ii) beforehand. Without doing this we can approximate them to only the order $O(j^2/T)$; this order is, however, enough to show the assertion (ii). More precisely, similar to (a) - (c) without $(*)$ we get
\begin{align*}
&\bE \brac{ \bra{ \int^{t_{i+j}}_{t_i} \lambda(t_i/T) dl}^2 - \bra{ \int^{t_{i+j}}_{t_i} \lambda(l/T) dl}^2  \con \cF_{t_i,T}} \\
&= \ \bE \Biggl[ \underbrace{ \int^{t_{i+j}}_{t_i}  \set{ \lambda(t_i/T) -\lambda(l/T)}dl}_{\text{see i) without $(*)$}} \cdot \underbrace{ \int^{t_{i+j}}_{t_i} \set{ \lambda(t_i/T) + \lambda(l/T)} dl}_{=O(t_{i+j} -t_i)}  \con \cF_{t_i,T} \Biggr] \\
&= \ O\bra{ \frac{j^{2+\gamma'}}{T^{\gamma'}}I_{\set{m'= 0}} + \frac{j^3}{T} I_{\set{m' \neq 0}} }
\end{align*}
by Lemma 1. In particular,
\begin{align*}
\bE \brac{\lambda^2(t_i/T) (t_{i+j} - t_i)^2 \con \cF_{t_i,T}} & = \ j^2 + j + O\bra{ \frac{j^{2+\gamma'}}{T^{\gamma'}}I_{\set{m'= 0}} + \frac{j^3}{T}I_{\set{m' \neq 0}}}. 
\end{align*} 
\qed 


\noindent \textbf{Proof of Theorem~\ref{Thm:TickTimeVolEstimator}.} 
This proof is analogous to that of Theorem~\ref{Thm:ClockTimeVolatilityEstimator} with an exception of the bias term $(IV)$, which is non-trivial in this case. Therefore, we have omitted its details and give only the outline of the proof. Let 
\begin{align*}
\widehat{\sigma}^2_{pre}(u_o) &:=  \frac{T}{N}\frac{1}{g_2} \sum^M_{j=-M} k\bra{  \frac{j}{M}} \bra{\ot{Y}_{t_{i_o + jH},T}}^2 \\
&\qquad - \frac{T}{2NH} \frac{\sum^{H-1}_{l=1} h^2(l/H)}{g_2} \sum^N_{i=-N} k\bra{\frac{i}{N}} \bra{Y_{t_{i_o+i},T} - Y_{t_{i_o+i-1},T}}^2,
\end{align*}
where $M=M(T)=\lfloor N/H \rfloor$ (we might assume that $M=N/H$ is an integer). Similar to the previous theorem, we point out that the difference between the two statistics - $\widehat{\sigma}^2_{pre}(u_o)$ and $\widehat{\sigma}^2_{pavg}(u_o)$ - is asymptotically negligible, i.e. 
$$
\sqrt{N/T^{1/2}} \set {\widehat{\sigma}^2_{pre}(u_o) - \widehat{\sigma}^2_{pavg}(u_o)} = o_p(1). 
$$
This means seeking the limit distribution of $\widehat{\sigma}^2_{pre}(u_o)$ is sufficient to infer the limit of $\widehat{\sigma}^2_{pavg}(u_o)$. Likewise, we must show the following statements:
\begin{align*}
&(I): \sqrt{\frac{N}{H}} \Biggl( \frac{T}{N}\frac{1}{g_2} \sum^M_{j=-M} k\bra{  \frac{j}{M}} \set{ \bra{\ot{X}_{t_{i_o + jH},T}}^2 - \bE \brac{ \bra{\ot{X}_{t_{i_o + jH},T}}^2 \con N_{\cdot,T}  } } \Biggr)  \xrightarrow{\cD} \ \cN(0, \ \xi^2_A), \\
&(II): \sqrt{\frac{NH}{T}} \Biggl( \frac{T}{N}\frac{2}{g_2} \sum^M_{j=-M} k\bra{ \frac{j}{M}}  \bra{\ot{X}_{t_{i_o + jH},T}} \bra{\ot{\ep}_{i_o + jH} } \Biggr)  \xrightarrow{\cD} \ \cN(0, \ \xi^2_B),\\
&(III): \sqrt{\frac{NH^3}{T^2}} \Biggl( \frac{T}{N}\frac{1}{g_2} \sum^M_{j=-M} k\bra{  \frac{j}{M}} \set{ \bra{\ot{\ep}_{i_o + jH}}^2 - \bE \bra{\ot{\ep}_{i_o + jH}}^2 }  \Biggr)  \xrightarrow{\cD} \ \cN(0, \ \xi^2_C),
\end{align*}
and the asymptotic biases
\begin{align*}
&(IV): \sqrt{\frac{N}{H}} \Biggl( \frac{T}{N}\frac{1}{g_2} \sum^M_{j=-M} k\bra{  \frac{j}{M}}   \bE \brac{ \bra{\ot{X}_{t_{i_o + jH},T}}^2 \con N_{\cdot,T} } -\sigma^2(u_o)  \Biggr) = \ o_p(1), \quad \text{ and } \\
&(V): \sqrt{\frac{NH^3}{T^2}} \Biggl( \frac{T}{N}\frac{1}{g_2} \sum^M_{j=-M} k\bra{  \frac{j}{M}} \bE \brac{ \bra{\ot{\ep}_{i_o + jH}}^2} \\
&\qquad \qquad -\frac{T}{2NH} \frac{\sum^{H-1}_{l=1} h^2(l/H)}{ g_2 } \sum^N_{i=-N} k\bra{\frac{i}{N}} \bra{Y_{t_{i_o+i},T} - Y_{t_{i_o+i-1},T}}^2   \Biggr) = \ o_p(1).
\end{align*}
As was pointed out, one of the most difficult parts in this proof is the derivation of $(IV)$, which is not obvious in the tick-time volatility estimation. More precisely, we have
\begin{align*}
&\sqrt{\frac{N}{H}} \set{ \frac{T}{N}\frac{1}{g_2} \sum^M_{j=-M} k\bra{  \frac{j}{M}}   \bE \brac{ \bra{\ot{X}_{t_{i_o + jH},T}}^2 \con N_{\cdot,T} } -\sigma^2(u_o)  } \\
&\quad = \sqrt{\frac{N}{H}} \frac{T}{N}\frac{1}{g_2} \frac{1}{H}\sum^N_{j=-N} \set{ k\bra{ \frac{\lfloor j/H \rfloor}{M} }-k\bra{  \frac{j}{N}} } \bE \brac{ \bra{\ot{X}_{t_{i_o + \ubra{j}},T}}^2 \con N_{\cdot,T} }  \\
&\qquad +\sqrt{\frac{N}{H}} \set{ \frac{T}{N}\frac{1}{g_2}\frac{1}{H}\sum^N_{j=-N} k\bra{  \frac{j}{N}}  \bE \brac{ \bra{\ot{X}_{t_{i_o + \ubra{j}},T}}^2 \con N_{\cdot,T} } \ -\sigma^2(u_o)  } \\
&\quad=:  (i)+(ii),
\end{align*}
where
\begin{align*}
& (i) \leq \ C \cdot \sqrt{\frac{N}{H}} \frac{T}{N} \frac{1}{H}\sum^N_{j=-N} k^{(1)}(...)\abs{\frac{\lfloor j/H \rfloor}{M} - \frac{j}{N}} \cdot\sum^{H-1}_{l=1 } g^2\bra{\frac{l}{H}} \frac{1}{T}  = \ O\bra{\sqrt{\frac{H}{N}}} = \ o(1)
\end{align*}
by the boundedness of derivatives $k'(\cdot)$, and
\begin{align*}
(ii)=& \sqrt{\frac{N}{H}} \Biggl( \frac{1}{NHg_2} \sum^N_{j=-N} k\bra{  \frac{j}{N}} \sum^{H-1}_{l=1} g^2\bra{\frac{l}{H}}  \sigma^2\bra{\frac{t_{i_o+\ubra{j}+l}}{T}} -\sigma^2(u_o) \Biggr) \\
=& \ O\bra{\sqrt{\frac{N}{H}} \cdot \abs{\frac{H}{T}}^\gamma} + \underbrace{ \sqrt{\frac{N}{H}}  \frac{1}{N} \sum^N_{j=-N} k\bra{  \frac{j}{N}} \set{  \sigma^2\bra{\frac{t_{i_o+j}}{T}} -\sigma^2(u_o) } }_{(\clubsuit)} = \ o(1)  +  o_p(1),
\end{align*}
as the segment condition holds. To obtain $(\clubsuit) = o_p(1)$, we need to apply Corollary~1 many times. Here we demonstrate $(\clubsuit) $ only for the case $m,m'=2$  to simplify notation (other cases are analogous). For $\sigma(\cdot) \in \cC^{2,\gamma}$ and $\lambda(\cdot) \in \cC^{2,\gamma'}$ we can expand
\begin{align*} 
&\sqrt{\frac{N}{H}}\bra{ \frac{1}{N} \sum^N_{j=-N} k\bra{\frac{j}{N}} \set{\sigma^2\bra{ \frac{t_{i_o-j}}{T}} -\sigma^2\bra{\frac{t_o}{T}} } } \notag\\
&= \sqrt{\frac{N}{H}} \frac{1}{N}  \sum^N_{j=-N} k\bra{\frac{j}{N}} \Biggl(   (\sigma^2(u_o))^{(1)}\bra{ \frac{t_{i_o-j}}{T} - \frac{t_o}{T}}  + \frac{ ( \sigma^2(u_o))^{(2)}}{2} \bra{ \frac{t_{i_o-j}}{T} - \frac{t_o}{T}}^2 + O \bra{ \abs{ \frac{t_{i_o-j}}{T} - \frac{t_o}{T}}^{2+\gamma} }  \Biggr) \\
&=:  (A) + (B) + (C). 
\end{align*}
It is easy to see that 
$$
(A) = \sqrt{\frac{N}{H}} \frac{1}{N} (\sigma^2(u_o))^{(1)}  \sum^N_{j=-N} k\bra{\frac{j}{N}}  \bra{ \frac{t_{i_o-j}}{T} - \frac{t_{i_o}}{T} } + o_p(1) =: (A_1) + o_p(1)
$$ 
We will show that $\bE\brac{ (A_1)^2} = o(1)$, which results that~$(A_1)$ is $o_p(1)$. By the symmetry of $k$ we can rewrite 
$$
(A_1) = \sqrt{\frac{N}{H}} \frac{1}{N} (\sigma^2(u_o))^{(1)}  \sum^N_{j=1} k\bra{\frac{j}{N}}  \set{ \bra{ \frac{t_{i_o-j}}{T} - \frac{t_{i_o}}{T} } + \bra{\frac{t_{i_o+j}}{T} - \frac{t_{i_o}}{T}} },  
$$
thus (setting $\sigma' := (\sigma^2(u_o))^{(1)}$, $\sigma'' := (\sigma^2(u_o))^{(2)}  $, $\lambda' := (\lambda(u_o))^{(1)}$ and $\lambda :=\lambda(u_o)$)
\begin{align*}
\bE \brac{ (A_1)^2 }  &= \frac{N}{H}\frac{1}{N^2}  \sigma'^2 \sum^N_{i,j=1}  k\bra{\frac{i}{N}} k\bra{\frac{j}{N}} \times\\
&\quad \times \Biggl( \bE\brac{ \bra{\frac{t_{i_o-i}}{T} - \frac{t_{i_o}}{T}} \bra{\frac{t_{i_o-j}}{T} - \frac{t_{i_o}}{T}} }  + \bE\brac{ \bra{\frac{t_{i_o-i}}{T} - \frac{t_{i_o}}{T}} \bra{\frac{t_{i_o+j}}{T} - \frac{t_{i_o}}{T}}   } \\
&\qquad \quad+ \bE\brac{ \bra{\frac{t_{i_o+i}}{T} - \frac{t_{i_o}}{T}} \bra{\frac{t_{i_o-j}}{T} - \frac{t_{i_o}}{T}}   } + \bE\brac{ \bra{\frac{t_{i_o+i}}{T} - \frac{t_{i_o}}{T}} \bra{\frac{t_{i_o+j}}{T} - \frac{t_{i_o}}{T}}   } \Biggr) \\ 
&=: (I_1) + (I_2) + (I_3) + (I_4).
\end{align*}
Since $N_{\cdot,T}$ has independent increments, the non-overlapping interarrival times are independent, particularly $\bE\brac{ (t_i - t_j) (t_k - t_l)} = \bE\brac{t_i - t_j} \bE\brac{t_k - t_l}$ for $l <k\leq j < i$. Therefore, by Corollary~1 
\begin{align*}
(I_2) &= \frac{N}{H}\frac{1}{N^2}  \sigma'^2 \sum^N_{i,j=1}  k\bra{\frac{i}{N}} k\bra{\frac{j}{N}} \bE\brac{ \frac{t_{i_o-i}}{T} - \frac{t_{i_o}}{T} } \bE\brac{ \frac{t_{i_o+j}}{T} - \frac{t_{i_o}}{T}  } \\
&=  \frac{N}{H}\frac{1}{N^2T^2} \sigma'^2 \sum^N_{i,j=1}  k\bra{\frac{i}{N}} k\bra{\frac{j}{N}}  \set{ \frac{-i}{\lambda} -\frac{1}{2}\frac{\lambda'}{\lambda^2} \frac{i^2}{T} + O\bra{ \frac{i^3}{T^2}}}  \set{ \frac{j}{\lambda} -\frac{1}{2}\frac{\lambda'}{\lambda^2} \frac{j^2}{T} + O\bra{ \frac{j^3}{T^2}}}\\
&=: (I_{2,1}) + ...+ (I_{2,9}) 
\end{align*}
with 
\begingroup
\allowdisplaybreaks
\begin{align*}
(I_{2,1}) =& - \frac{N^3}{HT^2} \cdot \frac{\sigma'^2 }{\lambda^2} \bra{\frac{1}{N} \sum^N_{i=1} k\bra{ \frac{i}{N}} \frac{i}{N}}^2, \\
(I_{2,2}) =& \frac{N^4}{HT^3} \cdot \frac{\sigma'^2}{2} \frac{\lambda'}{\lambda^3} \bra{\frac{1}{N} \sum^N_{i=1} k\bra{ \frac{i}{N}} \frac{i}{N}} \bra{\frac{1}{N} \sum^N_{j=1} k\bra{ \frac{j}{N}} \frac{j^2}{N^2}}, \quad (I_{2,3}) = O\bra{ \frac{N^5}{HT^4} }, \\
(I_{2,4}) =& -\frac{N^4}{HT^3} \cdot \frac{\sigma'^2}{2} \frac{\lambda'}{\lambda^3} \bra{\frac{1}{N} \sum^N_{j=1} k\bra{ \frac{j}{N}} \frac{j}{N}} \bra{\frac{1}{N} \sum^N_{i=1} k\bra{ \frac{i}{N}} \frac{i^2}{N^2}}, \tag{ $(I_{2,2})$ cancels out $(I_{2,4})$}\\ 
(I_{2,5}) =& \ O\bra{\frac{N^5}{HT^4} }, \ (I_{2,6}) = \ O\bra{ \frac{N^6}{HT^5}  }, \\
(I_{2,7}) =&  \ O\bra{\frac{N^5}{HT^4} }, \ (I_{2,8}) = \ O\bra{ \frac{N^6}{HT^5}  } \ \text{and } (I_{2,9}) = \ O\bra{ \frac{N^7}{HT^6}}.
\end{align*}
\endgroup
In particular, we obtain $ (I_2) = (I_{2,1}) + o(1)$ if $N^5/(HT^4) \to 0$ (which is satisfied by our segment conditions). We will see later that $(I_{2,1})$ is eliminated so that $(I_2) = o(1)$. Let us continue to $(I_4)$.
\begin{align*}
(I_4) =& \ \frac{N}{HT^2} \cdot \sigma'^2\frac{1}{N^2} \sum^N_{i=1}  k^2\bra{\frac{i}{N}} \bE\brac{ (t_{i_o+i} - t_{i_o})^2} \\
&+ \frac{N}{HT^2} \cdot \sigma'^2\frac{1}{N^2} \sum_{i\neq j, j < i} k\bra{\frac{i}{N}} k\bra{\frac{j}{N}}  \bE \brac{ (t_{i_o+i} -t_{i_o})(t_{i_o+j} - t_{i_o}) }\\
&+ \frac{N}{HT^2} \cdot \sigma'^2\frac{1}{N^2} \sum_{i\neq j, j > i} k\bra{\frac{i}{N}} k\bra{\frac{j}{N}}  \bE \brac{ (t_{i_o+i} -t_{i_o})(t_{i_o+j} - t_{i_o}) }\\
=:& \ (I_{4,1}) + (I_{4,2}) + (I_{4,3}),
\end{align*}
where
\begin{equation*}
(I_{4,1}) = \frac{N}{HT^2} \cdot \sigma'^2\frac{1}{N^2} \sum^N_{i=1}  k^2\bra{\frac{i}{N}} \set{ \frac{i^2+i}{\lambda^2} + O\bra{ \frac{i^3}{T}}} \to \ 0 \tag{by Corollary 1}
\end{equation*}
and
\begin{align*}
(I_{4,2}) &= \frac{N}{HT^2} \cdot \sigma'^2\frac{1}{N^2} \sum_{i\neq j, j < i} k\bra{\frac{i}{N}} k\bra{\frac{j}{N}} \bE \brac{ (t_{i_o+i} -t_{i_o+j})(t_{i_o+j} - t_{i_o}) + (t_{i_o+j} - t_{i_o})^2 } \\
 &=: (I_{4,2,1}) + (I_{4,2,2}).
\end{align*}
Again by the independence of interarrival times and Corollary~1 we get
\begin{align*}
(I_{4,2,1}) 
&= \frac{N}{HT^2}  \sigma'^2\frac{1}{N^2} \sum_{i\neq j, j < i} k\bra{\frac{i}{N}} k\bra{\frac{j}{N}} \set{  \bE \brac{ t_{i_o+i} -t_{i_o}} \bE\brac{t_{i_o+j} - t_{i_o} }  - \bra{\bE\brac{t_{i_o+j} - t_{i_o}} }^2  } \\
&= \frac{N^3}{HT^2}  \frac{\sigma'^2}{\lambda^2} \frac{1}{N^2} \sum_{i\neq j, j < i} k\bra{\frac{i}{N}} k\bra{\frac{j}{N}} \frac{i}{N} \frac{j}{N} - \frac{N^3}{T^2} \frac{\sigma'^2}{\lambda^2} \frac{1}{N^2} \sum_{i\neq j, j < i} k\bra{\frac{i}{N}} k\bra{\frac{j}{N}} \frac{j^2}{N} + o(1)\\
&=: (I_{4,2,1,1}) + (I_{4,2,1,2}) + o(1)
\end{align*}
as $N^4/(HT^3) \to 0$ (in the same way as in the derivation of $(I_2)$). The next term
\begin{align*}
(I_{4,2,2}) =& \ \frac{N}{HT^2} \cdot \sigma'^2\frac{1}{N^2} \sum_{i\neq j, j < i} k\bra{\frac{i}{N}} k\bra{\frac{j}{N}} \set{ \frac{j^2 +j}{\lambda^2_.} + O\bra{ \frac{j^3}{T}}} \\
=& \ \frac{N}{HT^2} \cdot \frac{\sigma'^2}{\lambda^2} \sum_{i\neq j, j < i} k\bra{\frac{i}{N}} k\bra{\frac{j}{N}} \frac{j^2}{N^2} + o(1)
\end{align*}
as $N^4/(HT^3) \to 0$. Evidently, the first term of $(I_{4,2,2})$ wipes out $(I_{4,2,1,2})$. Furthermore, by changing the role of $i$ and $j$ the same result can be derived for $(I_{4,3})$, therefore $(I_4) = (I_{4,2,1,1}) +(I_{4,3,1,1}) + o(1)$. Fortunately we see that the sum of $(I_{4,2,1,1})$  and $(I_{4,3,1,1})$ is exactly the negative of $(I_{2,1})$, thus
$$
(I_2) + (I_4) = o(1).
$$

Analogously, we can show that $(I_1) + (I_3) = o(1)$, so we can conclude that $(A1)$ is $o_p(1)$. Finally, the negligiblities of $(B)$ and $(C)$ is verified by $\bE \abs{(B)} = o(1)$ and $\bE \abs{(C)} = o(1)$. 
\qed\\


\noindent \textbf{Proof of Theorem~\ref{Thm:AlternativeClockTimeVolEstimator}.} 
\noindent In the case of $m,m'=0$, it is necessary to restrict the the bandwidth size $\fb$ and the segment length $N$ according to Theorems~\ref{Thm:IntensityEstimator} and \ref{Thm:TickTimeVolEstimator} respectively, i.e.
$$
N/T^{1/2}= o\bigl( T^{\frac{\gamma}{1+2\gamma}} \bigr) \quad \text{and} \quad \fb T= o\bigl( T^{\frac{2\gamma'}{1+2\gamma'}} \bigr)
$$
in order to obtain those limit distributions. We analyze it by looking at the following situations: \\

\noindent (a) If $\gamma'= \frac{\gamma}{2\gamma +2}$ then $N/T^{1/2}=O(\fb T) $. This gives 
\begin{align*}
&\sqrt{\fb T} \set{ \widehat{\sigma}^2_{pavg}(u_o) \widehat{\lambda}(u_o)  - \sigma^2(u_o)\lambda(u_o)} \\
&\quad = \underbrace{\widehat{\sigma}^2_{pavg}(u_o)}_{  \xrightarrow\bP \sigma^2(u_o)} \cdot \sqrt{\fb T} \set{ \widehat{\lambda}(u_o) - \lambda(u_o)}  +  \lambda(u_o) \underbrace{\frac{\sqrt{\fb T}}{\sqrt{N/T^{1/2}}}}_{= \sqrt{c_1}} \cdot \sqrt{N/T^{1/2}}\set{ \widehat{\sigma}^2_{pavg}(u_o)  - \sigma^2(u_o)}\\
&\quad \xrightarrow\cD \ \cN\bra{ 0,  \sigma^4(u_o)\lambda(u_o)\int_\bR \fK^2(x)dx  +  c_1 \lambda^2(u_o)\set{\delta\xi^2_A + \frac{1}{\delta}\xi^2_B + \frac{1}{\delta^3}\xi^2_C} }, 
\end{align*}
since the limits of the first and second terms are independent.\\

\noindent (b) If $\gamma'> \frac{\gamma}{2\gamma +2}$ then $N/T^{1/2} = o(\fb T)$. This leads to 
\begin{align*}
\sqrt{N/T^{1/2}} \set{ \widehat{\sigma}^2_{pavg}(u_o) \widehat{\lambda}(u_o)  - \sigma^2(u_o)\lambda(u_o)}  &= \ \widehat{\lambda}(u_o) \cdot \sqrt{N/T^{1/2}}\set{ \widehat{\sigma}^2_{pavg}(u_o) - \sigma^2(u_o)}  \\ 
&\xrightarrow\cD \cN\bra{ 0, \lambda^2(u_o)\set{\delta\xi^2_A  + \frac{1}{\delta}\xi^2_B + \frac{1}{\delta^3}\xi^2_C} }.
\end{align*}

\noindent (c) If $\gamma'< \frac{\gamma}{2\gamma +2}$ then $\fb T = o(N/T^{1/2})$, which implies that
\begin{align*}
&\sqrt{\fb T} \set{ \widehat{\sigma}^2_{pavg}(u_o) \widehat{\lambda}(u_o)  - \sigma^2(u_o)\lambda(u_o)} \xrightarrow\cD \ \cN\bra{ 0, \ \sigma^4(u_o)\lambda(u_o)\int_\bR \fK^2(x)dx}.
\end{align*}
Other cases ($m=0$ together with $m'=1,2$; or $m=1,2$) can be verified by the same arguments with the corresponding restrictions on $\fb$ and $N$.
\qed 

\normalsize

\section*{REFERENCES}

\small
\begin{description}
\baselineskip1.3em
\itemsep-0.04cm

\item A\"it-Sahalia, Y., and Jacod, J. (2014), {\it High-Frequency Financial Econometrics}, Princeton, New Jersey: Princeton University Press.
\item An\'e, T., and Geman, H. (2000),
	``Order Flow, Transaction Clock, and Normality of Asset Returns,"{\it The Journal of Finance}, 55, 2259--2284.
\item Bandi, F. M., and Russell, J. R. (2008),
	``Microstructure Noise, Realized Variance, and Optimal Sampling,"{\it The Review of Economic Studies}, 75, 339-369.
\item Barndorff-Nielsen, O. E., Hansen, P. R., Lunde, A., and Shephard, N. (2008),
	``Designing Realized Kernels to Measure the ex post Variation of Equity Prices in the Presence of Noise,"{\it Econometrica}, 76, 1481--1536.
\item Belomestny, D. (2011),
	``Statistical Inference for Time-Changed L\'evy Processes via Composite Characteristic Function Estimation,"{\it The Annals of Statistics}, 39, 2205--2242.	
\item Bibinger, M., Jirak, M., and Rei\ss, M. (2015),
	``Volatility Estimation Under One-sided Errors With Applications to Limit Order Books,"{\it ArXiv Preprint} [online], arXiv:1408.3768v4. Available at \url{http://arxiv.org/pdf/1408.3768v4.pdf}.
\item Br\'emaud, P. (1981), {\it Point Processes and Queues}, New York: Springer-Verlag.
\item Clark, P. K. (1973),
	``A Subordinated Stochastic Process Model With Finite Variance for Speculative Prices,"{\it Econometrica}, 41, 135--155.
\item Dahlhaus, R. (1997),
	``Fitting Time Series Models to Nonstationary Processes,"{ \it The Annals of Statistics}, 25, 1--37.
\item Dahlhaus, R., and Neddermeyer, J. C. (2013),
	``Online Spot Volatility-Estimation and Decomposition with Nonlinear Market Microstructure Noise Models,''{\it Journal of Financial Econometrics}, 12, 174--212.
\item Delbaen, F., and Schachermayer, W. (1994),
	``A General Version of the Fundamental Theorem of Asset Pricing,"{\it Mathematische Annalen}, 300, 463--520.	
\item Gabaix, X., Gopikrishnan, P., Plerou, V., and Stanley, H.E. (2003),
	``A Theory of Power-Law Distributions in Financial Market Fluctuations,"{\it Nature}, 423, 267--270.
\item Griffin, J. E., and Oomen, R. C. (2008),
	``Sampling Returns for Realized Variance Calculations: Tick Time or Transaction Time?,"{\it Econometric Reviews}, 27, 230--253.
\item Hall, P., and Heyde, C. C. (1980), {\it Martingale Limit theory and its Application}, New York: Academic Press.
\item Hansen, P. R., and Lunde, A. (2006),
	``Realized Variance and Market Microstructure Noise,"{\it Journal of Business} \&{\it Economic Statistics}, 24, 127--161.
\item Jacod, J., Li, Y., Mykland, P. A., Podolskij, M., and Vetter, M. (2009),
	``Microstructure Noise in the Continuous Case: The Pre-Averaging Approach,"{\it Stochastic Processes and their Applications}, 119, 2249--2276.
\item Jacod, J., and Protter, P. E. (2012),
	{\it Discretization of processes}. Berlin Heidelberg: Springer-Verlag.
\item Jones, C. M., Kaul, G., and Lipson, M. L. (1994),
	``Information, Trading, and Volatility,"{\it Journal of Financial Economics}, 36, 127--154.
\item Koo, B., and Linton, O. (2012),
	``Estimation of Semiparametric Locally Stationary Diffusion Models,"{\it  Journal of Econometrics}, 170, 210--233.
\item Kuo, H.-H. (2006), {\it Introduction to Stochastic Integration}, New York, Springer.
\item Munk, A., and Schmidt-Hieber, J. (2010),
	``Lower Bounds for Volatility Estimation in Microstructure Noise Models,"{\it Borrowing Strength: Theory Powering Applications -A Festschrift for Lawrence D. Brown}, 6, 43--55.
\item Plerou, V., Gopikrishnan, P., Gabaix, X., A Nunes Amaral, L., and Stanley, H.E. (2001),
   ``Price Fluctuations, Market Activity and Trading Volume,"{\it Quantitative Finance}, 1, 262--269.
\item Podolskij, M., and Vetter, M. (2009),
	``Estimation of Volatility Functionals in the Simultaneous Presence of Microstructure Noise and Jumps,"{\it Bernoulli}, 15, 634--658.
\item Roueff, F., von Sachs, R., and Sansonnet, L. (2016), ``Locally stationary Hawkes processes,"{\it   Stochastic Processes and their Applications}, 126, 1710--1743.
\item Reiss, M. (2011),
	``Asymptotic Equivalence for Inference on the Volatility From Noisy Observations,"{\it The Annals of Statistics}, 39, 772--802.
\item Zhang, L. (2006),
	``Efficient Estimation of Stochastic Volatility Using Noisy Observations: A Multi-Scale Approach,"{\it Bernoulli}, 12, 1019--1043.
\item Zhang, L., Mykland, P. A., and A\"it-Sahalia, Y. (2005),
	``A Tale of Two Time Scales: Determining Integrated Volatility With Noisy High-Frequency Data,"{\it Journal of the American Statistical Association}, 100, 1394--1411.

\end{description}

\end{document}